%% file: Gilpin_covprop_juq_paper.tex
\documentclass[review,onetabnum]{siamonline190516}
\usepackage{amsfonts}
\usepackage{algorithmic}
\usepackage{graphicx,epstopdf} 
\usepackage[caption=false]{subfig}

\input{Gilpin_covprop_juq_shared}

\ifpdf
\hypersetup{
  pdftitle={Continuum covariance propagation for understanding variance loss in advective systems},
  pdfauthor={S. Gilpin, T. Matsuo, and S. E. Cohn}
}
\fi

\newcommand{\bs}{\boldsymbol}

\begin{document}
\nolinenumbers

\maketitle

\begin{abstract}
Motivated by the spurious variance loss encountered during covariance propagation in atmospheric and other large-scale data assimilation systems, we consider the problem for state dynamics governed by the continuity and related hyperbolic partial differential equations. This loss of variance is often attributed to reduced-rank representations of the covariance matrix, as in ensemble methods for example, or else to the use of dissipative numerical methods. Through a combination of analytical work and numerical experiments, we demonstrate that significant variance loss, as well as gain, typically occurs during covariance propagation, even at full rank. The cause of this unusual behavior is a discontinuous change in the continuum covariance dynamics as correlation lengths become small, for instance in the vicinity of sharp gradients in the velocity field. This discontinuity in the covariance dynamics arises from hyperbolicity: the diagonal of the kernel of the covariance operator is a characteristic surface for advective dynamics. Our numerical experiments demonstrate that standard numerical methods for evolving the state are not adequate for propagating the covariance, because they do not capture the discontinuity in the continuum covariance dynamics as correlations lengths tend to zero. Our analytical and numerical results demonstrate in the context of mass conservation that this leads to significant, spurious variance loss in regions of mass convergence and gain in regions of mass divergence. The results suggest that developing local covariance propagation methods designed specifically to capture covariance evolution near the diagonal may prove a useful alternative to current methods of covariance propagation.

\end{abstract}



\section{Introduction}\label{sec:intro}

At the heart of modern data assimilation is covariance propagation. Data assimilation techniques evolve the estimation error covariance along with the model state, either explicitly as in the Kalman filter, implicitly as in variational methods, or using a reduced-rank approximation as in ensemble schemes \cite{Ka60, TaCo87, CoTa87, Ev94}. To provide context for the problem addressed in this paper, we start with a stochastic model state $N$-vector $\bs{q}$, which is propagated discretely in data assimilation schemes from time $t_{k-1}$ to $t_k$ as
\begin{equation}\label{eq:discrete state}
    \bs{q}_k = \bs{M}_{k,k-1}\bs{q}_{k-1},
\end{equation}
where $\bs{M}_{k,k-1}$ is the deterministic $N\times N$ propagation matrix representing the model dynamics. For simplicity, we consider here the linear case with no forcing, random or otherwise. From the model state we can define the $N\times N$ symmetric positive semi-definite covariance matrix at time $t_k$, 
\begin{equation}
    \bs{P}_k = \mathbb{E}[(\bs{q}_k-\overline{\bs{q}}_k)(\bs{q}_k-\overline{\bs{q}}_k)^T],
\end{equation} 
where $\mathbb{E}[\cdot]$ is the expectation operator,  $\overline{\bs{q}}_k = \mathbb{E}[\bs{q}_k]$ is the mean state, and superscript $T$ denotes transpose. The basic equation of discrete covariance propagation behind modern data assimilation schemes then follows directly from the discrete state propagation \cref{eq:discrete state},
\begin{equation}\label{eq:discrete covariance}
	\bs{P}_{k} = \bs{M}_{k,k-1}\bs{P}_{k-1}\bs{M}_{k,k-1}^T,
\end{equation}
where $\bs{P}_{k-1}$ and $\bs{P}_k$ are the covariance matrices at times $t_{k-1}$ and $t_{k}$, respectively \cite[chapter 6]{Ka60, Ja70}. We omit a process noise term in \cref{eq:discrete covariance}, which will be discussed later.

Motivated by atmospheric data assimilation schemes used for global numerical weather prediction \cite{Da93,Ka03}, we consider covariance propagation associated with hyperbolic partial differential equations (PDEs). Let $\bs{x} \in S^2_r$ where $S^2_r$ is the surface of the sphere of radius $r$ and take time $t\geq0$. To fix ideas, we focus on the continuity equation, which describes the continuum evolution of the stochastic model state $q = q(\bs{x},t)$ as follows:
\begin{gather}
	q_t + \bs{v}\bs{\cdot}\bs{\nabla}q + (\bs{\nabla \cdot v})q = 0, \nonumber \\
	q(\bs{x}, t_0) = q_0(\bs{x}). \label{eq:continuity S2}
\end{gather}
The subscript $t$ denotes the time derivative unless noted otherwise.
The two-dimensional velocity field $\bs{v} = \bs{v}(\bs{x},t)$ is taken to be deterministic and continuously differentiable while the initial state $q_0$ is stochastic with mean $\overline{q}_0$. \Cref{eq:continuity S2} is the statement of mass conservation for a passively advected tracer in a thin layer of atmosphere between two isentropic surfaces, for instance \cite[sections 2.5\,--\,2.6]{MeCoChLy00,HoHa13}. For the model state in \cref{eq:continuity S2}, we can define the covariance between two points $\bs{x}_1, \bs{x}_2 \in S^2_r$ as
\begin{equation}\label{eq:covariance def}
    P(\bs{x}_1,\bs{x}_2,t) = \mathbb{E}[(q(\bs{x}_1,t) - \overline{q}(\bs{x}_1,t))(q(\bs{x}_2,t) - \overline{q}(\bs{x}_2,t))].
\end{equation}
The corresponding covariance evolution equation for \cref{eq:continuity S2} on $S^2_r\times S^2_r$ and for $t \geq 0$ becomes
\begin{gather}
	P_t + \bs{v}_1\bs{\cdot\nabla}_1P + \bs{v}_2\bs{\cdot\nabla}_2P + (\bs{\nabla}_1\bs{\cdot v}_1 + \bs{\nabla}_2\bs{\cdot v}_2)P = 0, \nonumber\\
	P(\bs{x}_1,\bs{x}_2,t_0) = P_0(\bs{x}_1,\bs{x}_2). \label{eq:covariance S4} 
\end{gather}
Here $\bs{\nabla}_i$ denotes the gradient with respect to $\bs{x}_i$, and $\bs{v}_i = \bs{v}(\bs{x}_i,t)$, $i=1,2$.

Both the model state and covariance equations given in \cref{eq:continuity S2,eq:covariance S4} are hyperbolic PDEs. For the covariance equation, the hyperplane $\bs{x}_1=\bs{x}_2$ is everywhere characteristic, so that solutions on that hyperplane are independent of initial conditions away from that hyperplane \cite[p.\,3130]{Co93}. As a result, we show in \cref{sec:analysis} that there is a discontinuous change in solutions to \cref{eq:covariance S4} along the $\bs{x}_1,\bs{x}_2$-hyperplane in the limit as correlation lengths approach zero, for example in the vicinity of sharp gradients in the velocity field which can arise naturally, as seen for instance in \cite{LyCoZhChMeOlRe04}. 

First, suppose the initial covariance $P_0$ is continuous on $S_r^2\times S_r^2$, and denote it by $P_0^d(\bs{x}_1,\bs{x}_2)$. Thus, the stochastic initial state $q_0$ is spatially correlated, with variance $\sigma_0^2(\bs{x}) = P_0^d(\bs{x},\bs{x})$. It follows that the solution of the covariance evolution equation \cref{eq:covariance S4} along the $\bs{x}_1,\bs{x}_2$-hyperplane corresponds to the variance $\sigma^2(\bs{x},t)=P(\bs{x},\bs{x},t)$ and satisfies
\begin{gather}
	\sigma^2_t + \bs{v}\bs{\cdot\nabla}\sigma^2 + 2(\bs{\nabla\cdot v})\sigma^2 = 0,\nonumber\\
	\sigma^2(\bs{x},t_0) = \sigma^2_0(\bs{x}) = P_0^d(\bs{x},\bs{x}), \label{eq:variance S2}
\end{gather}
where $\sigma^2 = \sigma^2(\bs{x},t)$ for $\bs{x} \in S^2_r$ and $t \geq 0$, which can be derived either from \cref{eq:covariance S4} or directly from \cref{eq:continuity S2}. 


Now, suppose the initial state $q_0$ is spatially uncorrelated, and denote its covariance function by $P_0^c(\bs{x}_1)\delta(\bs{x}_1,\bs{x}_2)$ where $\delta$ is the Dirac delta and $P_0^c$ is continuous on $S_r^1$. As we show in \cref{sec:continuum}, covariances that are initially restricted to the $\bs{x}_1,\bs{x}_2$-hyperplane remain so for all time, as this hyperplane is a characteristic surface. For $t\geq 0$ the solution to \cref{eq:covariance S4} is then given by $P(\bs{x}_1,\bs{x}_2,t) = P^c(\bs{x}_1,t)\delta(\bs{x}_1,\bs{x}_2)$ where $P^c$ satisfies 
\begin{gather}
	P^c_t + \bs{v}\bs{\cdot\nabla}P^c + (\bs{\nabla\cdot v})P^c = 0, \nonumber\\  
	P^c(\bs{x},t_0) = P^c_0(\bs{x}). \label{eq:continuous spectrum S2}
\end{gather}
Thus, near zero correlation length, the behavior of solutions of the covariance evolution equation along the $\bs{x}_1,\bs{x}_2$-hyperplane changes abruptly from that of the variance equation \cref{eq:variance S2} to that of \cref{eq:continuous spectrum S2}, in the case of a nonzero divergent velocity field, $\bs{\nabla\cdot v} \neq 0$.

The characteristic behavior of the $\bs{x}_1,\bs{x}_2$-hyperplane in the continuum does not translate into discrete space for typical discretizations \cref{eq:discrete state} of \cref{eq:continuity S2}. Diagonal elements of $\bs{P}_k$ in \cref{eq:discrete covariance} depend on off-diagonal elements of $\bs{P}_{k-1}$, and a diagonal initial covariance in \cref{eq:discrete covariance} does not remain diagonal for all time. In this paper, we study the behavior of solutions of the continuum covariance evolution equation \cref{eq:covariance S4} near the $\bs{x}_1,\bs{x}_2$-hyperplane and contrast this with the behavior of discretizations \cref{eq:discrete covariance} near the diagonal, using a combination of analytical and numerical methods. We conduct numerical experiments using a one-dimensional version of \cref{eq:continuity S2} and study the covariance and variance propagation as a function of correlation length scales of the initial covariance. We find in some cases that the variance propagated numerically according to \cref{eq:discrete covariance} bears little resemblance to that of the continuum variance dynamics \cref{eq:variance S2}. In particular, variances propagated numerically according to \cref{eq:discrete covariance} tend to be better approximations of \cref{eq:continuous spectrum S2} than of \cref{eq:variance S2} for short initial correlation lengths, quite independently of any numerical dissipation or dispersion effects. This property manifests itself as a large, spurious loss of variance in regions of mass convergence, and a spurious variance gain in regions of mass divergence, as our analytical and numerical results illustrate. This behavior is the result of the discontinuous change in the continuum covariance dynamics in the limit as the correlation length tends to zero.



Spurious loss of variance is well known in the data assimilation literature. How covariances propagated through data assimilation schemes tend to underestimate the exact error covariances has long been noted \cite[pp.\,23\,--\,24]{Ma82}, but variance loss has been discussed primarily in the context of ensemble schemes \cite{Ma82, Le74, Ev94}, where spurious variance loss can be attributed to the use of reduced-rank covariance representations \cite{FuBe07}. Several methods have been developed to address variance loss to prevent filter divergence, such as covariance inflation \cite[chapter 9.2]{AnAn99, MiHo00, Ma82} and a scale-selective generalization of covariance inflation \cite[section 2.4.4]{Co10}. Loss of variance is sometimes addressed through an artificial model error or process noise term  added to the discrete covariance propagation in \cref{eq:discrete covariance} \cite[sections 8.8\,--\,8.9]{Ja70}. Accurately estimating an appropriate model error/process noise term is difficult because spurious variance loss can be due to several different known and unknown sources, though it has been shown that adding a model error term can help rectify the negative impact of variance loss, for instance, by increasing ensemble spread in the case of the ensemble schemes \cite{MiHo00, HoMi05}. Stochastic parameterization of subgrid scale physics also helps to increase ensemble spread to prevent filter divergence \cite[p.\,567]{Be17}. For the purpose of illuminating a root cause of variance loss, we will consider only the unforced covariance dynamics and omit a model error term, artificial or not, in this work. Specifically, the focus of this paper is on the spurious loss (and gain) of variance associated with the peculiar, discontinuous limiting behavior of solutions of the continuum covariance evolution equation \cref{eq:covariance S4}. This is a full-rank effect. We believe ours is the first work to identify and study this effect.


The layout of this paper is as follows. In \cref{sec:analysis}, we consider a slightly more generalized form of the continuity equation to study the continuum covariance propagation. In \cref{sec:preliminaries} we establish the generalized problem, defining the necessary operators and associated PDEs used for the analysis. This is followed by \cref{sec:continuum}, which derives the generalized version of \cref{eq:continuous spectrum S2} and discusses the discontinuous change in the continuum dynamics as initial correlation lengths approach zero. The analysis sections are followed by numerical experiments, where we illustrate spurious loss and gain of variance in full-rank covariance propagation through a simple one-dimensional example. \Cref{sec:set up} details the experimental setup of the one-dimensional problem, describing the numerical propagation methods, discretization schemes, and initial covariances. \Cref{sec:results} summarizes the results from these numerical experiments. \Cref{sec:conclusions and discussion} contains concluding remarks, followed by \cref{sec:appendix a,sec:appendix b} which contain additional derivations. 

\section{Analysis}\label{sec:analysis}
To gain insight into the spurious loss of variance associated with full-rank numerical covariance propagation, we first study covariance propagation in the continuum. We will consider the state and covariance equations as PDEs with solutions in the Hilbert space $L^2$, define the associated linear operators, and use tools from functional analysis to study these equations and operators. This continuum analysis is crucial for interpreting the results of the numerical experiments given later in \cref{sec:experiments}. 

\subsection{Preliminaries}\label{sec:preliminaries}
Let $\Omega = S_r^2$ and take $\bs{x} \in \Omega$ and $t\geq0$. We will consider the generalized advection equation for the model state $q = q(\bs{x},t)$,
\begin{gather}
	q_t + \bs{v}\bs{\cdot}\bs{\nabla}q + bq = 0, \nonumber \\
	q(\bs{x}, t_0) = q_0(\bs{x}). \label{eq:generalized advection}
\end{gather}
Here $b = b(\bs{x},t)$ is a scalar, noting that setting $b = \bs{\nabla \cdot v}$ yields the continuity equation \cref{eq:continuity S2}. From \cref{eq:generalized advection} we have
\begin{equation}\label{eq:unique soln}
    \frac{d}{dt}\int_\Omega q^2 d\bs{x} + \int_\Omega (2b - \bs{\nabla \cdot v})q^2d\bs{x} = 0,
\end{equation}
and we assume that $2b - \bs{\nabla\cdot v} \in L^\infty(\Omega)$ and $q_0 \in L^2(\Omega)$ so that \cref{eq:generalized advection} has a unique solution $q \in L^2(\Omega)$ for all time. We write this solution as 
\begin{equation}\label{eq:continuum state prop}
q(\bs{x},t) = (\bs{\mathcal{M}}_tq_0)(\bs{x}),
\end{equation}
where $\bs{\mathcal{M}}_t\colon L^2(\Omega)\mapsto L^2(\Omega)$ is the solution operator of \cref{eq:generalized advection},
\begin{equation}\label{eq:fundamental solution operator}
(\bs{\mathcal{M}}_tf)(\bs{x}) = \int_\Omega M(\bs{x},t;\bs{\xi})f(\bs{\xi})d\bs{\xi}.
\end{equation}
The subscript $t$ on all operators denoted using the calligraphy font style, as in $\bs{\mathcal{M}}_t$ for example, indicates the operator evaluated at time $t$, not the time derivative. The kernel of the operator $\bs{\mathcal{M}}_t$, $M = M(\bs{x},t;\bs{\xi})$, is the fundamental solution of \cref{eq:generalized advection},
\begin{gather}
M_t + \bs{v}\bs{\cdot\nabla}M + bM = 0, \nonumber \\
M(\bs{x},t_0;\bs{\xi}) = \delta(\bs{x},\bs{\xi}), \label{eq:fundamental solution}
\end{gather}
where the initial condition is the Dirac delta. Here, we simply view the Dirac delta as the kernel of the identity operator $\bs{\mathcal{I}}\colon L^2(\Omega)\mapsto L^2(\Omega)$,
\begin{equation}\label{eq:identity operator}
    (\bs{\mathcal{I}}f)(\bs{x}) = f(\bs{x}) = \int_\Omega \delta(\bs{x},\bs{\xi})f(\bs{\xi})d\bs{\xi}.
\end{equation}

\Cref{eq:continuum state prop} is analogous to the discrete state propagation computed in data assimilation schemes. We can propagate our discrete state $\bs{q}$ in \cref{eq:discrete state} from time $t_0$ to $t_{k}$,
\begin{equation}\label{eq:discrete state 2}
    \bs{q}_{k} = \underbrace{\bs{M}_{k,k-1}\bs{M}_{k-1,k-2}...\bs{M}_{2,1}\bs{M}_{1,0}}_{\bs{M}_{k,0}}\bs{q}_0,
\end{equation}
which is the discrete version of \cref{eq:continuum state prop} evaluated at time $t=t_k$; the operator $\bs{\mathcal{M}}_{t_k}$ is the continuum version of the propagation matrix $\bs{M}_{k,0}$.

With the model state now defined, we can derive the corresponding covariance evolution equation for $P = P(\bs{x}_1, \bs{x}_2, t)$ with $\bs{x}_1,\bs{x}_2 \in \Omega$ and $t\geq 0$ using the definition given in \cref{eq:covariance def},
\begin{gather}
	P_t + \bs{v}_1\bs{\cdot\nabla}_1P + \bs{v}_2\bs{\cdot\nabla}_2P + (b_1+b_2)P = 0, \nonumber \\
	P(\bs{x}_1,\bs{x}_2,t_0) = P_0(\bs{x}_1,\bs{x}_2), \label{eq:continuum covariance}
\end{gather}
where again, $\bs{\nabla}_i$ refers to the gradient with respect to $\bs{x}_i$, and $\bs{v}_i = \bs{v}(\bs{x}_i,t), \ b_i = b(\bs{x}_i,t)$ for $i=1,2$.

The solution of the covariance evolution equation \cref{eq:continuum covariance} can be expressed using the fundamental solution operator $\bs{\mathcal{M}}_t$ and its adjoint $\bs{\mathcal{M}}_t^{\bs{*}}$. The adjoint fundamental solution operator is defined using the inner product over the Hilbert space $L^2(\Omega)$,
\begin{equation}\label{eq:adjoint definition}
(\bs{\mathcal{M}}_t^{\bs{*}} f,g)_2 = (f,\bs{\mathcal{M}}_t g)_2 \quad \forall f, g \in L^2(\Omega).
\end{equation}
The adjoint operator, $\bs{\mathcal{M}}_t^{\bs{*}}\colon L^2(\Omega)\mapsto L^2(\Omega)$, can be expressed as an integral operator whose kernel $M^*$ is the solution to the adjoint final value problem associated with \cref{eq:generalized advection},
\begin{equation}\label{eq:adjoint fundamental solution operator}
    (\bs{\mathcal{M}}_t^{\bs{*}} f)(\bs{\xi}) = \int_\Omega M^*(\bs{\xi}; \bs{x},t)f(\bs{x})d\bs{x}.
\end{equation}
With respect to the method of characteristics, the fundamental solution operator $\bs{\mathcal{M}}_t$ propagates the solution forward along the characteristics determined by departure points, and the adjoint fundamental solution operator $\bs{\mathcal{M}}_t^{\bs{*}}$ propagates the solution backwards along the characteristics determined by arrival points. This yields the symmetry property \cite[p.\,729]{CoHi62} that at any fixed time $t$ the kernels satisfy
\begin{equation}\label{eq:adjoint property}
M(\bs{x},t;\bs{\xi}) = M^*(\bs{\xi};\bs{x},t).
\end{equation}
Using \cref{eq:adjoint property}, we can express the covariance in terms of the kernels of the fundamental solution and adjoint fundamental solution operators,
\begin{equation}\label{eq:covariance integral solution}
P(\bs{x}_1,\bs{x}_2,t) = \int_\Omega \int_\Omega M(\bs{x}_1,t;\bs{\xi}_1)P_0(\bs{\xi}_1,\bs{\xi}_2)M^*(\bs{\xi}_2;\bs{x}_2,t)d\bs{\xi}_2d\bs{\xi}_1,
\end{equation}
or simply
\begin{equation}\label{eq:operator cov prop}
    \bs{\mathcal{P}}_t = \bs{\mathcal{M}}_t\bs{\mathcal{P}}_0\bs{\mathcal{M}}_t^{\bs{*}},
\end{equation}
where $\bs{\mathcal{P}}_0\colon L^2(\Omega)\mapsto L^2(\Omega)$ is the operator whose kernel is $P_0$,
\begin{equation}\label{eq:initial cov operator}
    (\bs{\mathcal{P}}_0f)(\bs{x}_1) = \int_\Omega P_0(\bs{x}_1,\bs{x}_2)f(\bs{x}_2) d\bs{x}_2,
\end{equation}
and $\bs{\mathcal{P}}_t\colon L^2(\Omega)\mapsto L^2(\Omega)$ is the resulting operator at time $t$,
\begin{equation}\label{eq:covariance operator}
    (\bs{\mathcal{P}}_tf)(\bs{x}_1) = \int_\Omega P(\bs{x}_1,\bs{x}_2, t)f(\bs{x}_2) d\bs{x}_2.
\end{equation}
Thus, the covariance evolution equation \cref{eq:continuum covariance} is interpreted as the evolution equation for the kernel of the covariance operator $\bs{\mathcal{P}}_t$. As with the state propagation, \cref{eq:operator cov prop} evaluated at time $t_k$ is the continuum version of the discrete covariance propagation,
\begin{equation}\label{eq:discrete covariance 2}
    \bs{P}_k = \bs{M}_{k,0}\bs{P}_0\bs{M}_{k,0}^T,
\end{equation}
following from \cref{eq:discrete covariance,eq:discrete state 2}.

We next define the (left) polar decomposition of the fundamental solution operator $\bs{\mathcal{M}}_t$, which will bring to light important properties of the covariance evolution that will be discussed in \cref{sec:continuum}. The polar decomposition is a canonical form for all bounded linear operators on Hilbert spaces \cite[p.\,197]{ReSi72}. It is the unique decomposition $\bs{\mathcal{M}}_t = \bs{\mathcal{D}}_t\bs{\mathcal{U}}_t$ where $\bs{\mathcal{D}}_t = (\bs{\mathcal{M}}_t\bs{\mathcal{M}}_t^{\bs{*}})^{1/2}$ and $\bs{\mathcal{U}}_t$ is a partial isometry. 

To derive the polar decomposition, we decompose the fundamental solution into 
\begin{equation}\label{eq:M decomp}
M(\bs{x},t;\bs{\xi}) = d(\bs{x},t)u(\bs{x},t;\bs{\xi})
\end{equation}
where $d = d(\bs{x},t)$ and $u = u(\bs{x},t;\bs{\xi})$ satisfy the following PDEs,
\begin{gather}
d_t+ \bs{v}\bs{\cdot\nabla}d + \left(b -\frac12 \bs{\nabla\cdot v}\right)d = 0, \nonumber \\
d(\bs{x},t_0) = 1, \label{eq:d equation}
\end{gather}
\begin{gather}
u_t + \bs{v}\bs{\cdot\nabla}u + \frac12(\bs{\nabla\cdot v})u = 0, \nonumber \\
u(\bs{x},t_0;\bs{\xi}) = \delta(\bs{x},\bs{\xi}). \label{eq:u equation}
\end{gather}

The solution $u$ of \cref{eq:u equation} is quadratically conservative, 
\begin{equation}\label{eq:energy conservation}
\frac{d}{dt}\int_\Omega u^2(\bs{x},t;\bs{\xi})d\bs{x} = 0,
\end{equation}
and therefore defines a bounded linear operator  $\bs{\mathcal{U}}_t\colon L^2(\Omega)\mapsto L^2(\Omega)$ whose kernel is the solution to \cref{eq:u equation},
\begin{equation}\label{eq:u operator}
(\bs{\mathcal{U}}_tf)(\bs{x}) = \int_\Omega u(\bs{x},t;\bs{\xi})f(\bs{\xi})d\bs{\xi}.
\end{equation}
The operator $\bs{\mathcal{U}}_t$ is an invertible isometry and therefore unitary.


We show in \cref{sec:appendix a} that the solution of \cref{eq:continuum covariance} with the initial condition $P_0(\bs{x}_1,\bs{x}_2) = \delta(\bs{x}_1,\bs{x}_2)$ is $P(\bs{x}_1,\bs{x}_2,t) = d^2(\bs{x}_1,t)\delta(\bs{x}_1,\bs{x}_2)$, where $d$ is the solution to \cref{eq:d equation}. In other words, according to \cref{eq:operator cov prop},  $\bs{\mathcal{M}}_t\bs{\mathcal{M}}_t^{\bs{*}}$ is in fact a multiplication operator. A multiplication operator $\bs{\mathcal{K}}\colon L^2(\Omega) \mapsto L^2(\Omega)$ is defined as
\begin{equation}\label{eq:multiplciation operator}
    (\bs{\mathcal{K}}f)(\bs{x}) = k(\bs{x})f(\bs{x}) = \int_\Omega k(\bs{\xi})\delta(\bs{x},\bs{\xi})f(\bs{\xi})d\bs{\xi},
\end{equation}
where $k(\bs{x})\in L^\infty(\Omega)$ is the multiplication function. The operator $\bs{\mathcal{M}}_t\bs{\mathcal{M}}_t^{\bs{*}}\colon L^2(\Omega)\mapsto L^2(\Omega)$,
\begin{equation}\label{eq:fundmental product operator}
(\bs{\mathcal{M}}_t\bs{\mathcal{M}}_t^{\bs{*}}f)(\bs{x}) = d^2(\bs{x},t)f(\bs{x}),
\end{equation}
is a multiplication operator, where the multiplication function satisfies the following differential equation,
\begin{gather}
d^2_t + \bs{v}\bs{\cdot\nabla}d^2 + \left(2b -\bs{\nabla\cdot v}\right)d^2 = 0, \nonumber \\
d^2(\bs{x},t_0) = 1. \label{eq:d2 equation}
\end{gather}

As the operator $\bs{\mathcal{M}}_t\bs{\mathcal{M}}_t^{\bs{*}}$ is non-negative, its square root exists, and we define the operator $\bs{\mathcal{D}}_t\colon L^2(\Omega)\mapsto L^2(\Omega)$ as this square root,
\begin{equation}\label{eq:d operator}
(\bs{\mathcal{D}}_tf)(\bs{x}) = ((\bs{\mathcal{M}}_t\bs{\mathcal{M}}_t^{\bs{*}})^{1/2}f)(\bs{x}) = d(\bs{x},t)f(\bs{x}),
\end{equation}
with the multiplication function for $\bs{\mathcal{D}}_t$ being the solution to \cref{eq:d equation}. Note that because the operator $\bs{\mathcal{D}}_t$ is a multiplication operator with a real-valued multiplication function, it is self-adjoint.

The decomposition \cref{eq:M decomp} of the kernel of the fundamental solution operator gives us the (left) polar decomposition of $\bs{\mathcal{M}}_t$,
\begin{equation}\label{eq:m polar decomp}
\bs{\mathcal{M}}_t = \bs{\mathcal{D}}_t\bs{\mathcal{U}}_t,
\end{equation}
with $\bs{\mathcal{D}}_t$ and $\bs{\mathcal{U}}_t$ defined above. In the next section, we use this polar decomposition to study the continuum covariance propagation. 

\subsection{Continuum covariance propagation}\label{sec:continuum}
The solution of the continuum covariance evolution equation \cref{eq:continuum covariance} depends on the initial covariance $P_0(\bs{x}_1,\bs{x}_2)$, and we show in this section how the behavior of the solution changes as the initial correlation length tends to zero. To do so, we will interpret the covariance function $P(\bs{x}_1,\bs{x}_2,t)$ as the kernel of the operator $\bs{\mathcal{P}}_t$ defined in \cref{eq:covariance operator}, whose evolution can be written in terms of the fundamental solution operator, its adjoint, and the initial covariance, as given in \cref{eq:operator cov prop}. 

We consider two cases for the initial covariance $P_0(\bs{x}_1,\bs{x}_2)$. First, assume the initial covariance is continuous on $\Omega\times\Omega$, and denote it as $P_0^d(\bs{x}_1,\bs{x}_2)$. Since $P_0^d$ is continuous, the solution to the covariance equation \cref{eq:continuum covariance} is a strong solution and has a bounded $L^2$-norm. Therefore the corresponding covariance operator, which we will denote as $\bs{\mathcal{P}}_t^d$, is a self-adjoint Hilbert-Schmidt operator \cite[pp.\,210\,--\,211]{ReSi72}. Hilbert-Schmidt operators are a subclass of the compact operators, and it follows from spectral theory that self-adjoint Hilbert-Schmidt operators only contain eigenvalues in their spectrum \cite[pp.\,230\,--\,232]{HuNa01} with the possible exception of zero in the continuous spectrum. The set of eigenvalues is often referred to as the discrete spectrum, hence we can refer to the covariance operator $\bs{\mathcal{P}}_t^d$ as the discrete spectrum covariance operator.

Now, consider the case where the initial state $q_0$ is spatially uncorrelated, whose covariance we represent by $P_0^c(\bs{x}_1)\delta(\bs{x}_1,\bs{x}_2)$, and assume that the function $P_0^c$ is continuous on $\Omega$. The Dirac delta in the initial covariance reduces the initial covariance operator \cref{eq:initial cov operator} to a multiplication operator; denote this operator as $\bs{\mathcal{P}}_0^c$. We can see how this impacts the corresponding covariance operator, which we denote as $\bs{\mathcal{P}}^c_t$, by applying the polar decomposition \cref{eq:m polar decomp} to \cref{eq:operator cov prop} with $\bs{\mathcal{P}}_0 = \bs{\mathcal{P}}_0^c$,
\begin{equation}\label{eq:cts spectrum cov operator deriv}
    \bs{\mathcal{P}}_t^c = \bs{\mathcal{D}}_t\bs{\mathcal{U}}_t\bs{\mathcal{P}}_0^c\bs{\mathcal{U}_t}^{\bs{*}}\bs{\mathcal{D}}_t = \bs{\mathcal{D}}_t\tilde{\bs{\mathcal{P}}_t}\bs{\mathcal{D}}_t,
\end{equation}
where 
\begin{equation}\label{eq:P tilde}
    \tilde{\bs{\mathcal{P}}_t} = \bs{\mathcal{U}}_t\bs{\mathcal{P}}_0^c\bs{\mathcal{U}}_t^{\bs{*}}.
\end{equation}
From the definition of $\bs{\mathcal{U}}_t$ in \cref{eq:u operator}, it follows that the operator $\tilde{\bs{\mathcal{P}}_t}\colon L^2(\Omega)\mapsto L^2(\Omega)$ has kernel $\tilde{P} = \tilde{P}(\bs{x}_1,\bs{x}_2,t)$ given by
\begin{gather}
    \tilde{P}_t + \bs{v}_1\bs{\cdot\nabla}_1\tilde{P} + \bs{v}_2\bs{\cdot\nabla}_2\tilde{P} + \frac{1}{2}(\bs{\nabla}_1\bs{\cdot v}_1 + \bs{\nabla}_2\bs{\cdot v}_2)\tilde{P} = 0, \nonumber \\
    \tilde{P}(\bs{x}_1,\bs{x}_2,t_0) = P_0^c(\bs{x}_1)\delta(\bs{x}_1,\bs{x}_2), \label{eq:P tilde equation}
\end{gather}
whose solution is shown in \cref{sec:appendix a} to be
\begin{equation}\label{eq:P tilde solution}
    \tilde{P}(\bs{x}_1,\bs{x}_2,t) = \tilde{P}^c(\bs{x}_1,t)\delta(\bs{x}_1,\bs{x}_2)
\end{equation}
where
\begin{gather}
    \tilde{P}^c_t + \bs{v\cdot\nabla}\tilde{P}^c = 0, \nonumber \\
    \tilde{P}^c(\bs{x},t_0) = P_0^c(\bs{x}). \label{eq:P tilde c}
\end{gather}
Thus, $\tilde{\bs{\mathcal{P}}_t}$ is a multiplication operator, and it follows from \cref{eq:d equation,eq:P tilde c} that $\bs{\mathcal{P}}_t^c$ in \cref{eq:cts spectrum cov operator deriv} is also a multiplication operator, 
\begin{equation}\label{eq:cts spectrum cov operator}
    (\bs{\mathcal{P}}_t^cf)(\bs{x}_1) = P^c(\bs{x}_1,t)f(\bs{x}_1) = \int_\Omega P^c(\bs{x}_1,t)\delta(\bs{x}_1,\bs{x}_2)f(\bs{x}_2) d\bs{x}_2,
\end{equation}
where $P^c = d^2\tilde{P}$ is the solution to the continuous spectrum equation 
\begin{gather}
	P^c_t + \bs{v}\bs{\cdot\nabla}P^c + (2b - \bs{\nabla\cdot v})P^c = 0, \nonumber \\  
	P^c(\bs{x},t_0) = P^c_0(\bs{x}). \label{eq:continuous spectrum}
\end{gather}
Since multiplication operators contain only a continuous spectrum \cite[pp.\,219, 240]{HuNa01}, we will refer to $\bs{\mathcal{P}}_t^c$ as the continuous spectrum covariance operator and to $P^c$ as the continuous spectrum solution.



Thus, we have shown that the solution of the covariance evolution equation \cref{eq:continuum covariance} for initial condition $P^c_0(\bs{x}_1)\delta(\bs{x}_1,\bs{x}_2)$ is $P(\bs{x}_1,\bs{x}_2,t) = P^c(\bs{x}_1,t)\delta(\bs{x}_1,\bs{x}_2)$; white noise evolved under the state dynamics \cref{eq:generalized advection} remains white. Further, we can see explicitly that the diagonal dynamics of the covariance are governed by the continuous spectrum equation \cref{eq:continuous spectrum} for spatially uncorrelated initial states, rather than by the variance equation
\begin{gather}
	\sigma^2_t + \bs{v}\bs{\cdot\nabla}\sigma^2 + 2b\sigma^2 = 0, \nonumber \\
	\sigma^2(\bs{x},t_0) = \sigma^2_0(\bs{x}) = P_0^d(\bs{x},\bs{x}), \label{eq:variance}
\end{gather}
which follows directly from \cref{eq:continuum covariance} when the initial covariance is continuous on $\Omega\times\Omega$ . Only in cases where the velocity field is divergence-free are the continuous spectrum equation and variance equation identical. The diagonal dynamics of the covariance are governed by the variance equation \cref{eq:variance} for all continuous initial covariances, independently of nonzero initial correlation lengths, but at zero correlation length the diagonal dynamics change abruptly to those of the continuous spectrum equation \cref{eq:continuous spectrum}. 

Covariances associated with spatially uncorrelated initial states, which correspond to covariances with zero initial correlation length scales, are limiting cases in our analysis, as well as in practice. Atmospheric wind fields, for example, have sharp vertical correlation structures relative to long horizontal correlations that need to be represented in covariances for numerical weather prediction models \cite{Ph86}, and wind shear leads to tracer correlations that shrink in the direction perpendicular to the flow such that they are no longer spatially resolved \cite{LyCoZhChMeOlRe04}. Through careful analysis and use of the polar decomposition, we are able to derive the discontinuous change in the diagonal dynamics as the correlation length approaches zero, which is not readily apparent when first considering the model state and covariance equations or when only considering these equations in discrete space.

\section{Numerical experiments}\label{sec:experiments}
The purpose of the numerical experiments is to examine the evolution of the variance during discrete covariance propagation and its relation to the discontinuous change in diagonal dynamics as the initial correlation length approaches zero we have demonstrated in the continuum.  To clearly illustrate the problems associated with full-rank covariance propagation, a simple example is used so that our numerical results can be compared to a known exact solution. We conduct these experiments for various types of initial covariances specified with different correlation kernel functions with varying correlation length scales. For each type, cases with spatially stationary and nonstationary variance are considered.

\subsection{Experimental setup}\label{sec:set up}
For these experiments, we will consider the one-dimensional version of the continuity equation \cref{eq:continuity S2} over the unit circle $S^1_1$. We take the advection speed to be independent of time and spatially-varying,
\begin{equation}\label{eq:1d advection speed}
    v(x) = \sin(x) +2.
\end{equation}
The solution to the one-dimensional continuity equation with this advection speed can be obtained explicitly using the method of characteristics and is used as reference for our experiments; see \cref{sec:appendix b} for further discussion.

The spatial domain $S^1_1$ is discretized on a uniform grid, $x_i = i\Delta x, \ i = 0, 1, ..., N-1$ where $N=200$ and $\Delta x = \frac{2\pi}{N}$. The time discretization is given by $t_k = k\Delta t$ where the time step $\Delta t$ is determined from the Courant-Friedrichs-Lewy number $\lambda$,
\begin{equation}\label{eq:cfl}
    \lambda = \max_{x\in[0,2\pi]}v(x)\frac{\Delta t}{\Delta x} = 3\frac{\Delta t}{\Delta x} \leq 1.
\end{equation}
For these experiments, we take $\lambda = 1$. We ran experiments with several other values of $\lambda < 1$ (not shown) and found it did not have a significant impact on the results.

\subsubsection{Numerical covariance propagation}\label{sec:numerical cov prop}
We use two methods of propagation to illustrate the impact of numerical schemes on the discrete covariance propagation \cref{eq:discrete covariance} and to leverage insights from continuum covariance analysis in the form of the polar decomposition. The covariance matrix is propagated discretely using either of two methods:
\begin{enumerate}
    \item Traditional propagation: the covariance is propagated as in \cref{eq:discrete covariance}, where the matrix $\bs{M}_{k,k-1}$ is the finite difference discretization of the fundamental solution equation \cref{eq:fundamental solution}.
    \item Polar decomposition propagation: the polar decomposition of the fundamental solution operator \cref{eq:m polar decomp} is discretized and used in place of the matrix $\bs{M}_{k,k-1}$ in \cref{eq:discrete covariance}.
\end{enumerate}

In the polar decomposition propagation, \cref{eq:m polar decomp} is discretized as follows. The operator $\bs{\mathcal{D}}_t$ of \cref{eq:d operator} is a self-adjoint multiplication operator, therefore when $t= t_k$ its corresponding discretization is the diagonal matrix $\bs{D}_{k,0}$ whose diagonal elements are $(\bs{D}_{k,0})_{ii} = d(x_i,t_k)$. Here $d(x_i,t_k)$ is the solution to \cref{eq:d equation} evaluated on the discrete spatial grid up to time $t_k$, which we generate using the exact solution to \cref{eq:d equation}, as discussed in \cref{sec:appendix b}. The discretization of the operator $\bs{\mathcal{U}}_t$ of \cref{eq:u operator} is done via finite differences to generate the propagation matrix $\bs{U}_{k,k-1}$ corresponding to \cref{eq:u equation}. The matrix $\bs{U}_{k,k-1}$, which propagates the solution from time $t_{k-1}$ to $t_k$, is independent of time by virtue of \cref{eq:1d advection speed}, therefore we will denote $\bs{U}_{k,k-1} =\bs{U}$ for simplicity. To compute the covariance matrix at time $t_k$, $\bs{P}_k$, we do not construct the matrix $\bs{M}_{k,k-1}$ explicitly using the polar decomposition, but instead compute the covariance as follows:
\begin{align}
    \bs{P}_k &= \bs{D}_{k,0}\bs{U}_{k,k-1}\bs{U}_{k-1,k-2}...\bs{U}_{1,0}\bs{P}_0\bs{U}_{1,0}^T...\bs{U}_{k,k-1}^T\bs{D}_{k,0}\nonumber \\
    &= \bs{D}_{k,0}\bs{U}^k\bs{P}_0(\bs{U}^T)^k\bs{D}_{k,0}. \label{eq:discrete cov pd prop}
\end{align}


The Lax-Wendroff \cite{LaWe60} and Crank-Nicolson \cite{CrNi47} finite difference schemes are used to generate the matrices $\bs{U}_{k,k-1}$ and $\bs{M}_{k,k-1}$ corresponding to their respective PDEs \cref{eq:u equation} and \cref{eq:fundamental solution}. Like $\bs{U}_{k,k-1}$, $\bs{M}_{k,k-1}$ is independent of time, therefore through the rest of the paper we will denote $\bs{M}_{k,k-1}$ simply as $\bs{M}$. We choose these two simple finite difference schemes to illustrate their contrasting behaviors, particularly when generating the matrix $\bs{U}$. According to the continuum analysis, the operator $\bs{\mathcal{U}}_t$ is unitary. The Crank-Nicolson discretization preserves the unitary property of $\bs{\mathcal{U}}_t$  because the scheme is quadratically conservative, but the Lax-Wendroff scheme does not. In the numerical results, \cref{sec:results}, propagated covariances are labeled by the finite difference scheme used to construct the matrices $\bs{M}$ and $\bs{U}$ (Lax-Wendroff or Crank-Nicolson) followed by the method of propagation (traditional or polar decomposition) as described in the beginning of this section.

\subsubsection{Initial Covariances}\label{sec:initial cov}
\begin{figure}[htbp]
    \centering
    \includegraphics[width=\linewidth]{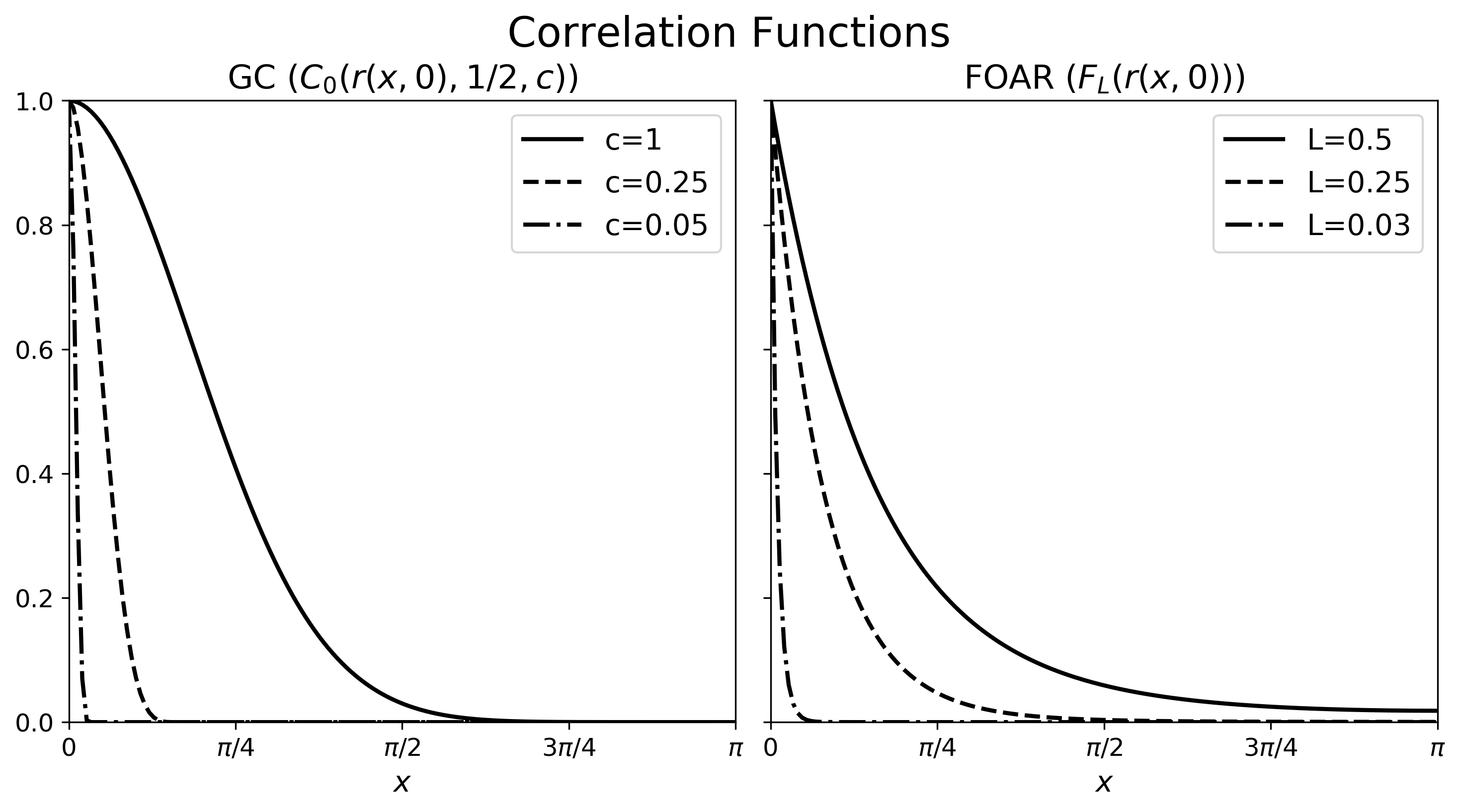}
    \caption{Examples of the correlation functions used to generate the initial covariances for numerical propagation. The GC correlation functions (left) are functions of compact support length parameter $c$, where the FOAR correlation functions (right) are functions of length scale $L$.}
    \label{fig:ex corr}
\end{figure}

We generate four types of initial covariances from the following two different correlation kernel functions: the Gaspari and Cohn (GC) $5^{th}$ order piecewise rational function \cite[Eq.~4.10]{GaCo99}, and the first order autoregressive function (FOAR) \cite[Eq.~(23)]{GaCoGuPa06}. Using each correlation function, we construct spatially stationary and nonstationary initial covariance matrices. The stationary initial covariances have an initial variance of one, while the spatially nonstationary initial covariances take the initial variance with spatially-varying standard deviation given as
\begin{equation}\label{eq:standard devation}
    \sigma_0(x) = \frac{\sin(3x)}{3} + 1.
\end{equation}
Nonstationarity in the initial covariances is introduced only through the variance. 

The GC correlation function $C_0(r(x_i,x_j),1/2,c)$ is a compactly supported approximation to a Gaussian function, supported on the interval $0 \leq r(x_i,x_j)\leq 2c$, where
\begin{equation}\label{eq:chordal distance}
    r(x_i,x_j) = 2\sin(|x_i-x_j|/2)
\end{equation} 
is the chordal distance between $x_i$ and $x_j$ on $S_1^1$. On the spatial grid for these experiments, values of $c=1$, 0.25, and 0.05 correspond to 100, 16, and 3 grid lengths ($\Delta x$), respectively, from the peak of the correlation function to where it becomes zero, and 33, 8, and just 1 grid length, respectively, from the peak value of 1 to values less than $0.2$; see \cref{fig:ex corr} for these examples.

The FOAR correlation function given by
\begin{equation}\label{eq:foar}
    F_L(r(x_i,x_j)) = exp(-r(x_i,x_j)/L)
\end{equation}
is continuous but non-differentiable at the origin because of its cusp-like behavior (see \cref{fig:ex corr}). As with the GC correlation function, the chordal distance \cref{eq:chordal distance} is used to reflect periodicity of the domain. The FOAR correlation functions are nonzero on the full spatial domain. On the spatial grid, $L=0.5$,  0.25, and 0.03 correspond to 26, 12, and just 1 grid length, respectively, from the peak of the correlation to where it becomes less than $0.2$.

\subsection{Experimental results}\label{sec:results}

\begin{figure}[htbp]
    \centering
    \includegraphics[width=\linewidth]{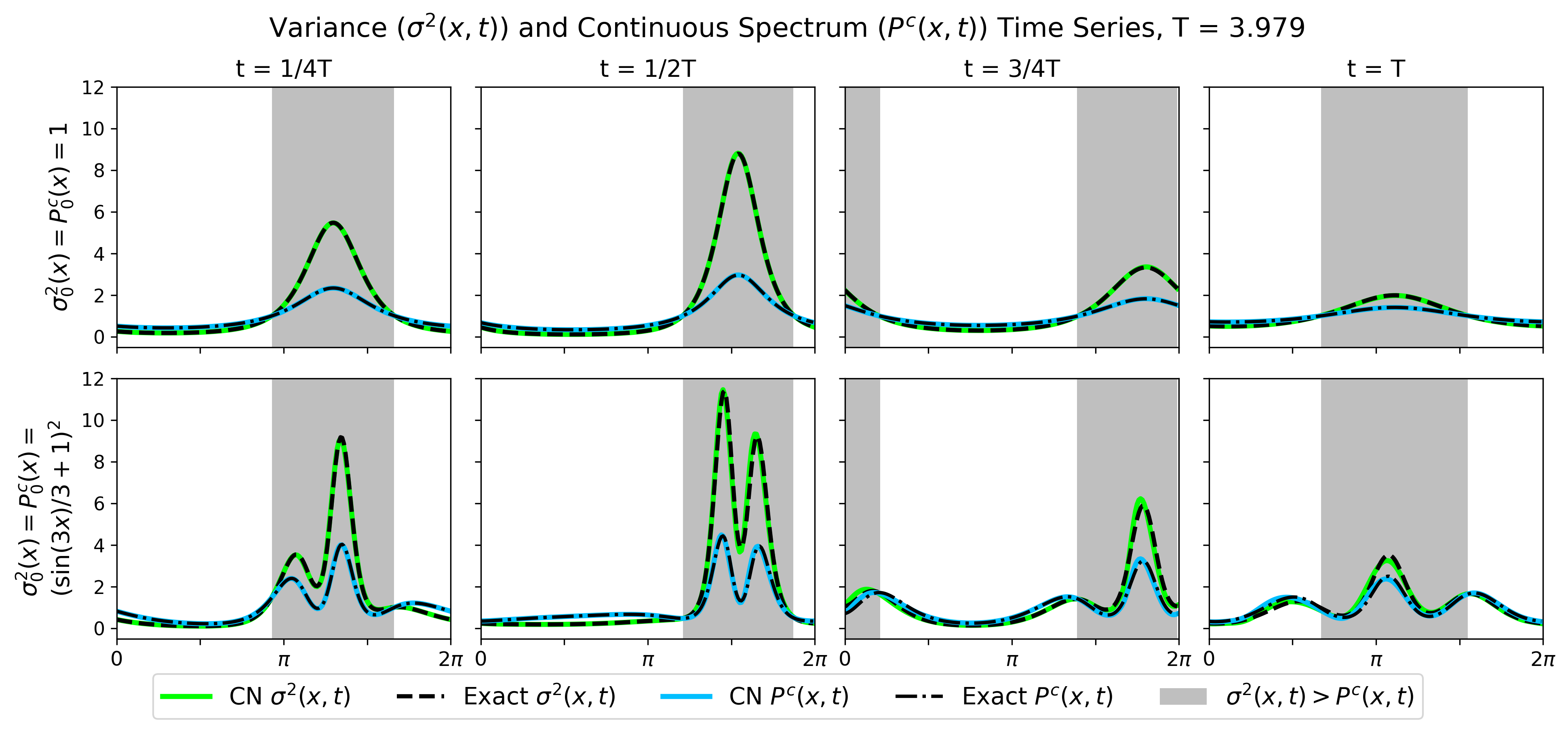}
    \caption{Top: Solutions to the one-dimensional variance equation \cref{eq:variance S2} and continuous spectrum equation \cref{eq:continuous spectrum S2} at various times with advection speed \cref{eq:1d advection speed}, for unit initial condition (top row) and spatially-varying initial condition (bottom row) taken to be the square of \cref{eq:standard devation}. Time $T=3.979$ corresponds to slightly after a full period. Exact solutions for each case are given in black dashed (variance) and black dot-dashed lines  (continuous spectrum). Green curves denote the solution to the variance equation computed with Crank-Nicolson (abbreviated CN); blue curves are solutions to the continuous spectrum equation computed with Crank-Nicolson. Regions highlighted in grey correspond to regions where $\sigma^2>P^c$; unhighlighted regions correspond to regions where $\sigma^2<P^c$.}
    \label{fig:variance and cts spec}
\end{figure}

From the continuum analysis, we expect the diagonal of the covariance matrix $\bs{P}_k$ to behave according to the variance equation for covariances with nonzero initial correlation lengths and according to the continuous spectrum equation for covariances with zero initial correlation length. To establish a baseline, \cref{fig:variance and cts spec} illustrates the solutions to the variance equation (1D version of \cref{eq:variance S2}) and continuous spectrum equation (1D version of \cref{eq:continuous spectrum S2}) for unit initial condition and for the spatially-varying initial condition given by the square of \cref{eq:standard devation}. Considering the exact solutions first (black), we see that the dynamics of the variance solution and continuous spectrum solution are quite different due to the spatially-varying advection speed \cref{eq:1d advection speed}. We also see that solving either the variance equation or continuous spectrum equation directly using the Crank-Nicolson scheme produces solutions that are nearly indistinguishable from the exact solutions. Solutions to the variance and continuous spectrum equations computed using Lax-Wendroff differ very slightly from Crank-Nicolson and are not shown.  

\begin{figure}[htbp]
    \centering
    \label{fig:GC1 a}\includegraphics[width=\linewidth]{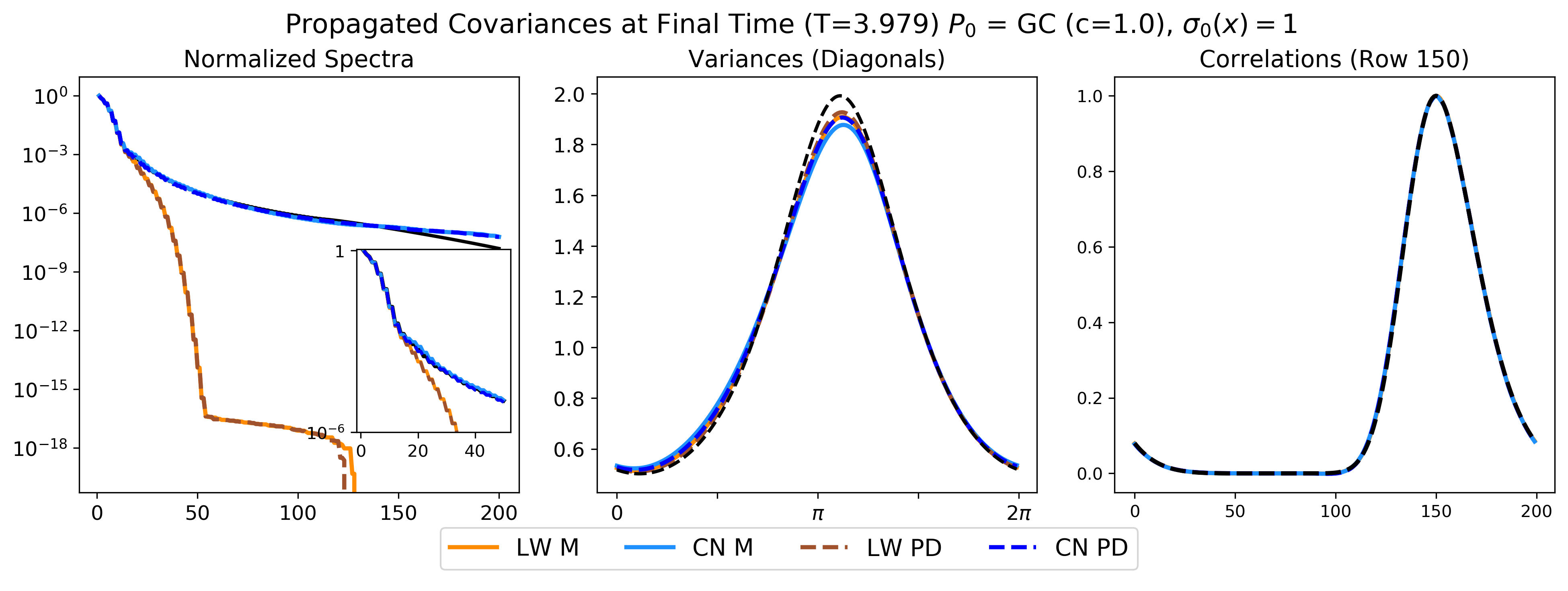}
    \label{fig:GC1 b}\includegraphics[width=\linewidth]{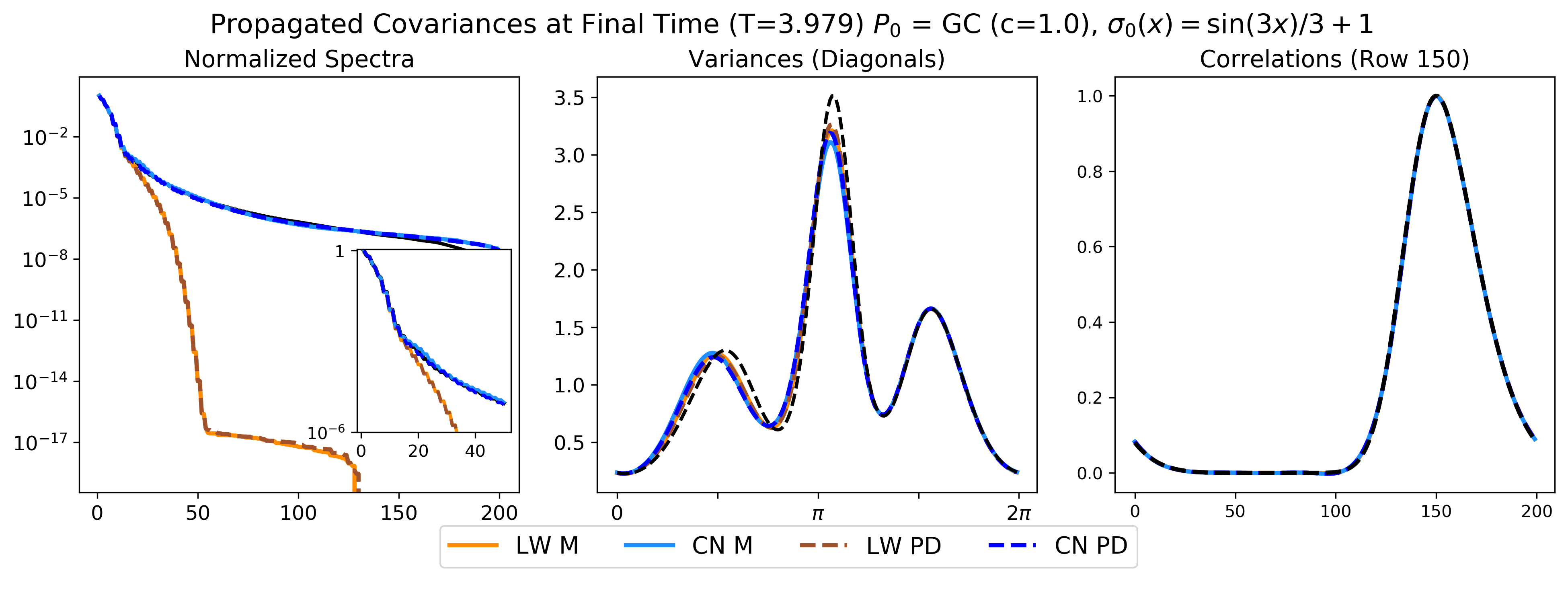}
    \caption{Propagated covariances at the final time $T$ (slightly after a full period) for GC initial correlations with $c=1$ (supported on the full domain) for all four propagation methods. Top row: unit initial variance. Bottom row: spatially-varying initial variance. The left panels correspond to the normalized spectra (relative to largest eigenvalue), the middle panels show the variances (diagonals), and the right panels show the correlations at row 150 (location of maximum in variance in space and time). The exact normalized spectra are given in solid black, and the exact variances and correlations are given in black dashed. Crank-Nicolson (CN) and Lax-Wendroff (LW) M refer to traditional propagation using Crank-Nicolson or Lax-Wendroff (solid). CN and LW PD refer to propagation using the polar decomposition (dashed) with the matrix  $\bs{U}$ constructed via the Crank-Nicolson or Lax-Wendroff scheme, respectively.}
    \label{fig:GC1}
\end{figure}
\Cref{fig:variance and cts spec} shows that propagating the diagonal of the covariance matrix $\bs{P}_k$ numerically, independently of the rest of the matrix, using the known diagonal dynamics of either the variance equation or the continuous spectrum equation, produces minimal discretization error. \Cref{fig:GC1} illustrates further that, at least for covariance matrices with relatively long initial correlation lengths, the numerically propagated full covariance itself also contains only minor discretization errors typically expected from finite difference approximations. \Cref{fig:GC1} shows results for the GC initial correlation supported on the full spatial domain ($c=1$). The dissipative behavior of both Lax-Wendroff schemes are clear in the normalized spectra (left panels), whereas both the traditional Crank-Nicolson propagation and polar decomposition propagation using the Crank-Nicolson $\bs{U}$ (hereafter referred to as Crank-Nicolson polar decomposition) capture the normalized spectra quite well. The small amount of variance loss and gain seen in both the constant initial variance and spatially-varying initial variance cases (middle panels of \cref{fig:GC1}) are consistent with dissipation and phase errors expected from finite differences. The polar decomposition methods (dashed) reduce the errors in the variance only slightly compared to the traditional methods (solid). The correlations at row 150 (right panels), which corresponds to where the variance reaches a maximum as the advection speed is at a minimum, are nearly identical to the exact correlations, suggesting minimal errors in correlation propagation.

\begin{figure}[htbp]
    \centering
    \includegraphics[width=\linewidth]{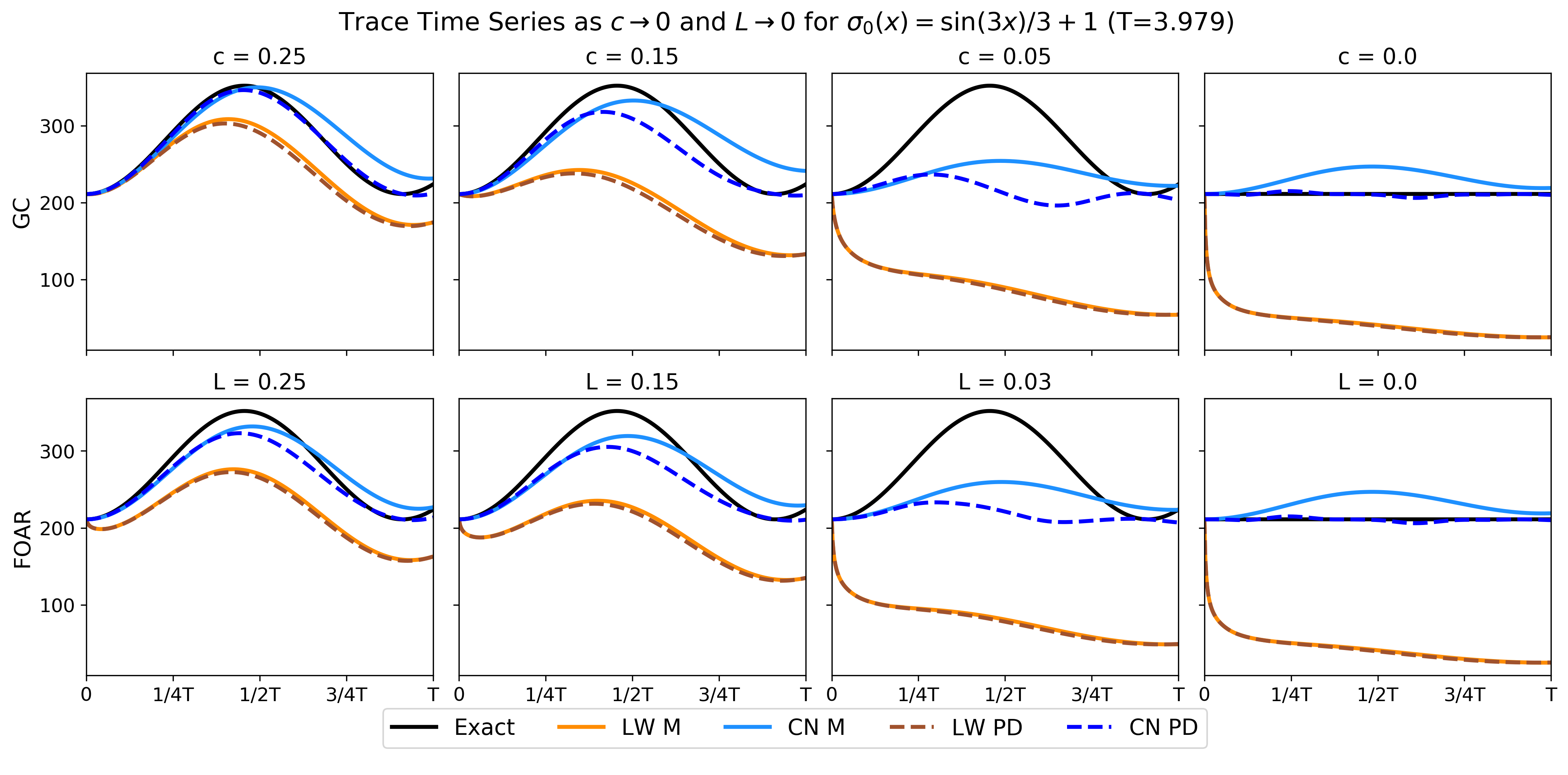}
    \caption{Trace time series for GC (top row) and FOAR (bottom row) for the spatially-varying initial variance taken to be the square of \cref{eq:standard devation}. Each panel corresponds to different values of $c$ and $L$, decreasing from left to right towards $c=L=0$. Refer to \cref{fig:GC1} for description of the curves. For the cases when $c=L=0$ (rightmost panels), the exact solutions (solid black) are constant in time due to the fact that continuous spectrum solution $P^c$ satisfies the continuity equation \cref{eq:continuous spectrum S2}, hence its integral over space is constant.}
    \label{fig:trace trig}
\end{figure}

The accuracy of numerically propagated full covariances is much worse for short initial correlation lengths, increasingly so as they approach zero. We monitor the total amount of variance lost or gained over time through the trace of the covariance matrix, $Tr(\bs{P}_k)$. \Cref{fig:trace trig} shows the trace time series for both the GC and FOAR cases with spatially-varying initial variance as initial correlation lengths tend to zero (the constant initial variance case results are similar to \cref{fig:trace trig} and not shown). The GC and FOAR cases in \cref{fig:trace trig} exhibit similar behaviors in the trace over time, even though their initial correlation structures are quite different. As the initial correlation lengths tend to zero, the amount of variance lost during propagation increases strikingly in the Lax-Wendroff schemes. The polar decomposition propagation schemes (dashed) are an improvement over traditional propagation in some cases but are worse in others, and generally suffer similar amounts of variance loss. We also observe that as initial correlation lengths tend to zero, the numerical schemes gradually approach their own limiting behavior at $c=L=0$, rather than a discontinuous change in dynamics as expected from the continuum analysis.

Had we not performed the continuum analysis of \cref{sec:analysis}, one might assume that the variance loss in \cref{fig:trace trig} is caused simply by dissipation. However, the Crank-Nicolson scheme is not dissipative and yet produces significant variance loss. In fact, we see for the traditional Crank-Nicolson propagation (solid light blue) that there are regions of both variance loss and gain. We also see that for short, nonzero initial correlation lengths ($c=0.05$, $L=0.03$ in particular), the numerical schemes better approximate the limiting case of $c=L=0$ than the correct behavior for $c,\ L >0$. This suggests that inaccurate variance propagation is particularly pronounced for short correlation lengths. 


The behavior of the trace time series in \cref{fig:trace trig} indicates that covariance propagation itself can be a source of spurious loss and gain of variance, however it does not indicate where exactly this manifests itself. To gain a better understanding of the source of variance loss and gain, we examine various aspects of the propagated covariance matrix for different types of initial covariances as we did in \cref{fig:GC1}, but now for covariances with shorter initial correlation lengths. 

\begin{figure}[htbp]
    \centering
    \label{fig:GC0.25 a}\includegraphics[width=\linewidth]{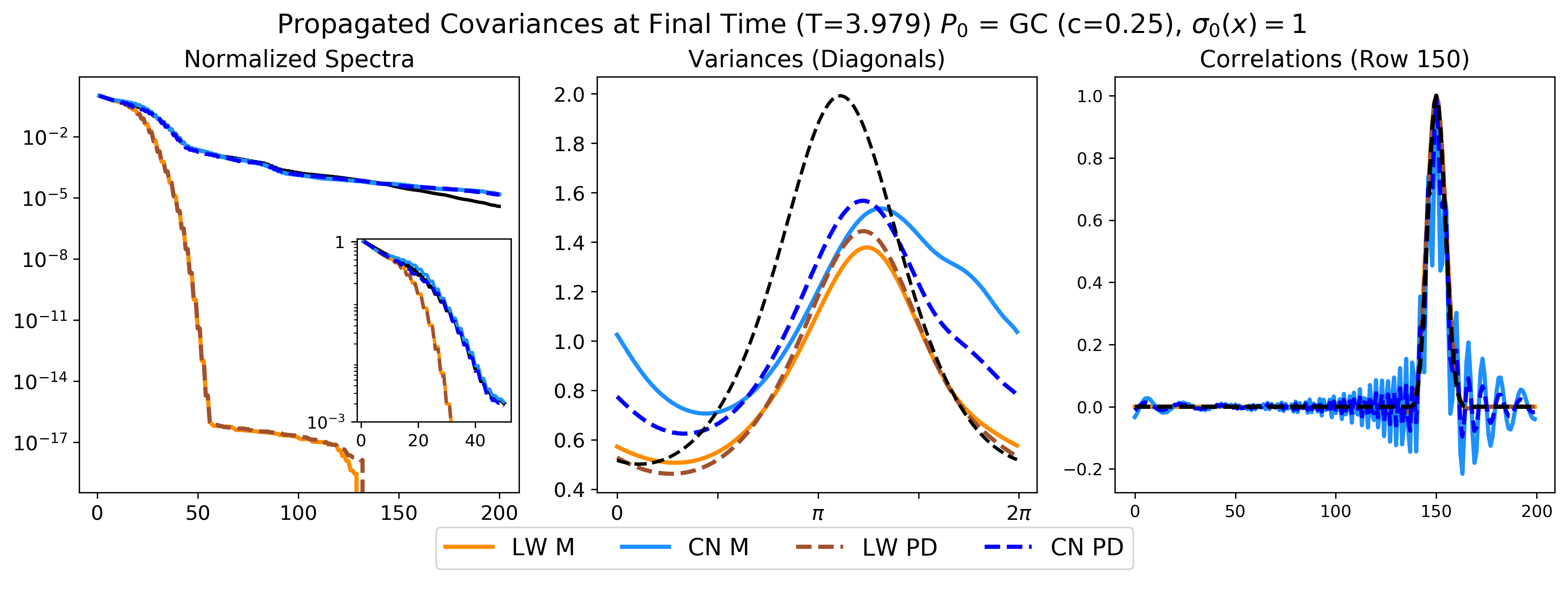}
    \label{fig:GC0.25 b}\includegraphics[width=\linewidth]{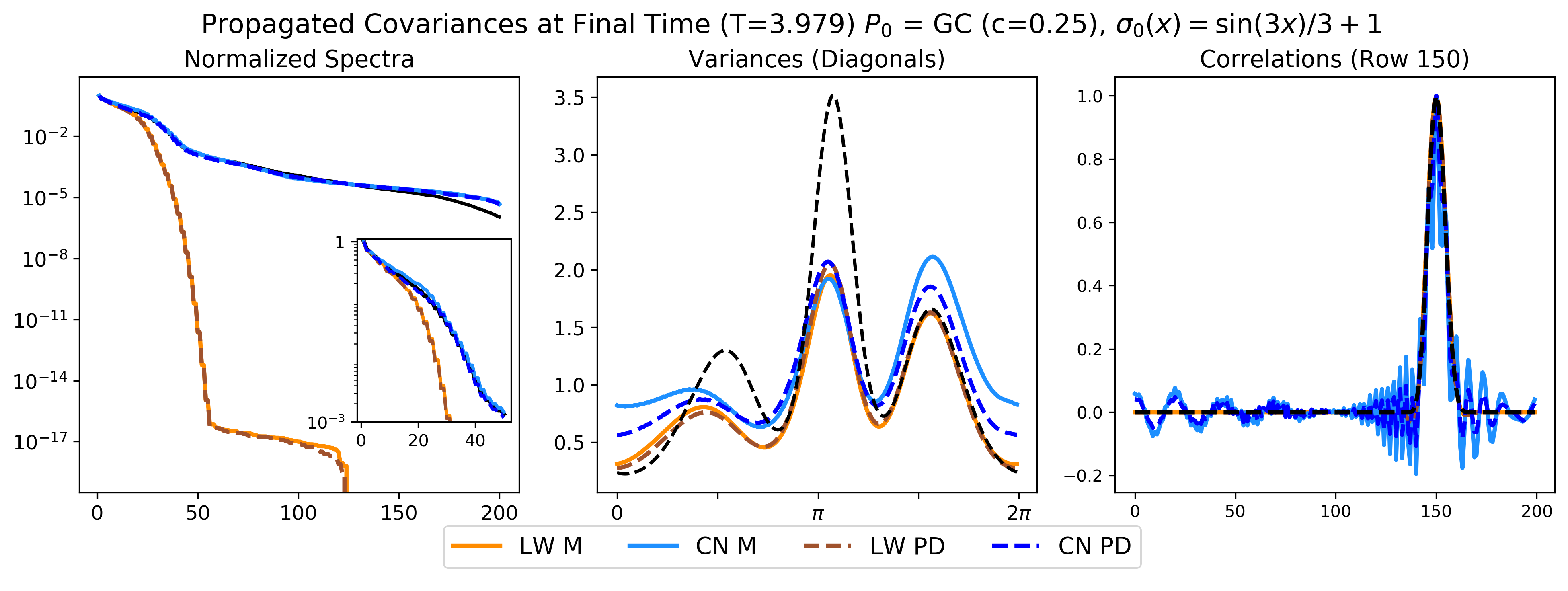}
    \caption{Same as \cref{fig:GC1} for $c=0.25$.}
    \label{fig:GC0.25}
\end{figure}

\begin{figure}[htbp]
    \centering
    \label{fig:FOAR0.25 a}\includegraphics[width=\linewidth]{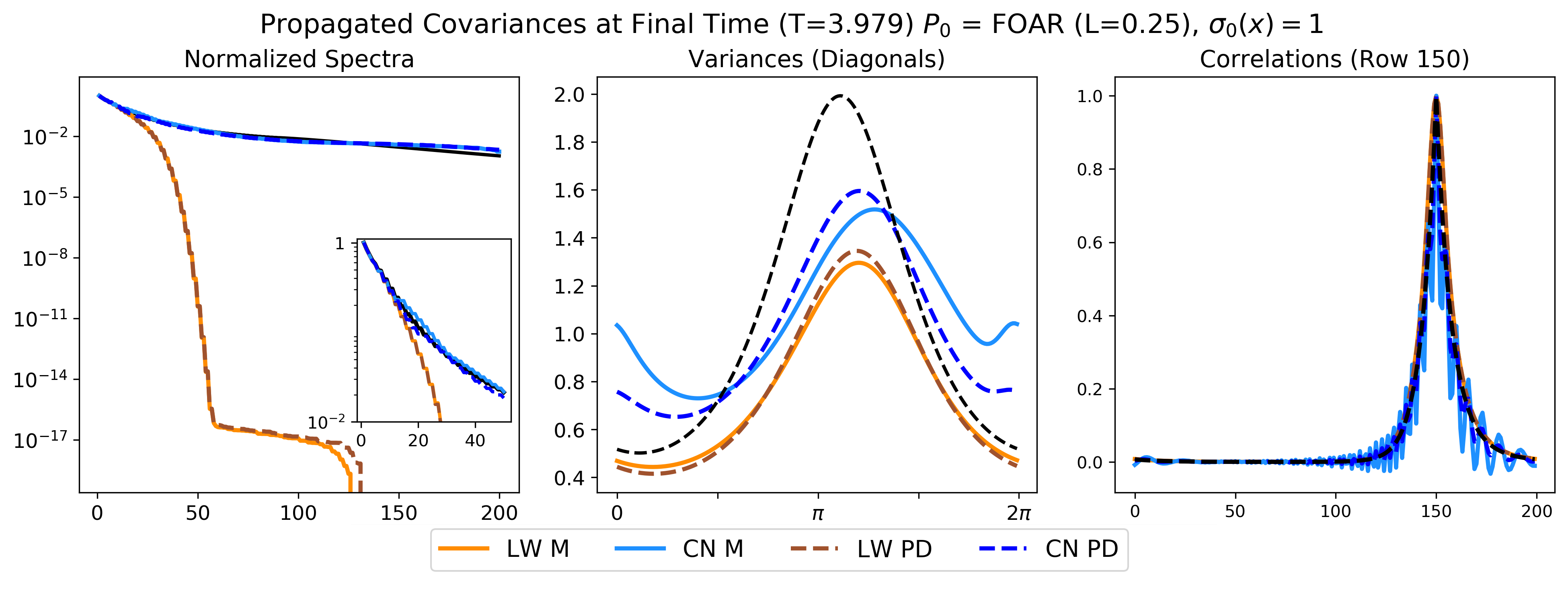}
    \label{fig:FOAR0.25 b}\includegraphics[width=\linewidth]{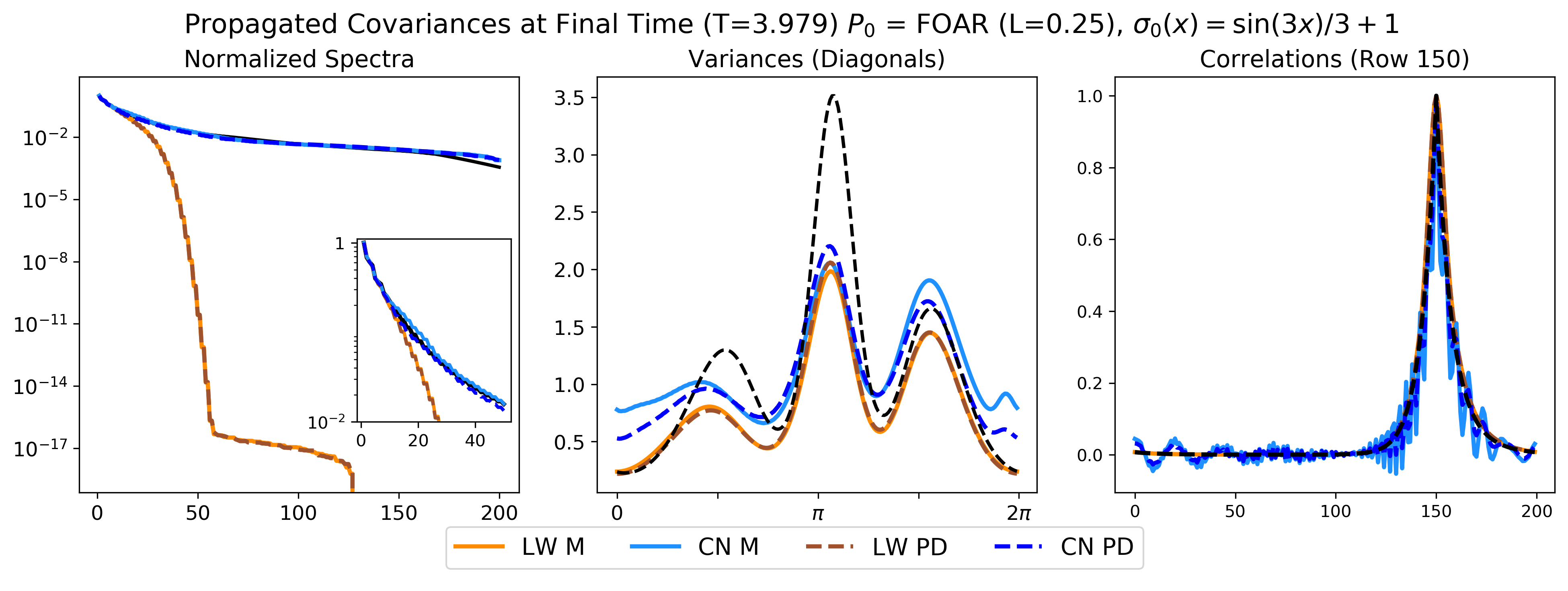}
    \caption{Same as \cref{fig:GC0.25} for the FOAR correlation function with $L=0.25$.}
    \label{fig:FOAR0.25}
\end{figure}

\Cref{fig:GC0.25,fig:FOAR0.25} are final time snapshots (in the same format as shown in \cref{fig:GC1}) of propagated covariances specified using the GC correlation function with $c=0.25$ and FOAR correlation function with $L=0.25$, respectively, which are initially well-resolved as described in \cref{sec:initial cov} and correspond to the mildest variance loss cases shown in leftmost panels of \cref{fig:trace trig}. The normalized spectra in both of these cases are similar to the spectra seen in \cref{fig:GC1}, where the Lax-Wendroff schemes are severely dissipative and both Crank-Nicolson methods approximate the exact spectrum moderately well. 

The variances extracted from the numerically propagated covariances in \cref{fig:GC0.25,fig:FOAR0.25} are strikingly different from the exact solution. Though covariances with these values of $c$ and $L$ are well-resolved initially, the resulting variance curves are smooth but wholly inaccurate. Across all four methods of propagation we see regions of both variance loss and variance gain, clearly illustrating that inaccurate variance propagation is the problem more so than dissipation. The variance propagation becomes worse when the initial covariance has a spatially-varying variance (bottom rows of \cref{fig:GC0.25,fig:FOAR0.25}), which is a more realistic situation in practice. The errors we observe in the variances, interestingly, are not reflected in the normalized spectra; considering the normalized spectra alone would not even hint at the problems occurring in the variance propagation. The correlations, as expected, show dispersion off the diagonal. The oscillatory behavior of both Crank-Nicolson schemes due to numerical dispersion is expected for this finite difference scheme \cite[p.\,46]{BeKe97}.


\begin{figure}[htbp]
    \centering
    \label{fig:tri ID a}\includegraphics[width=\linewidth]{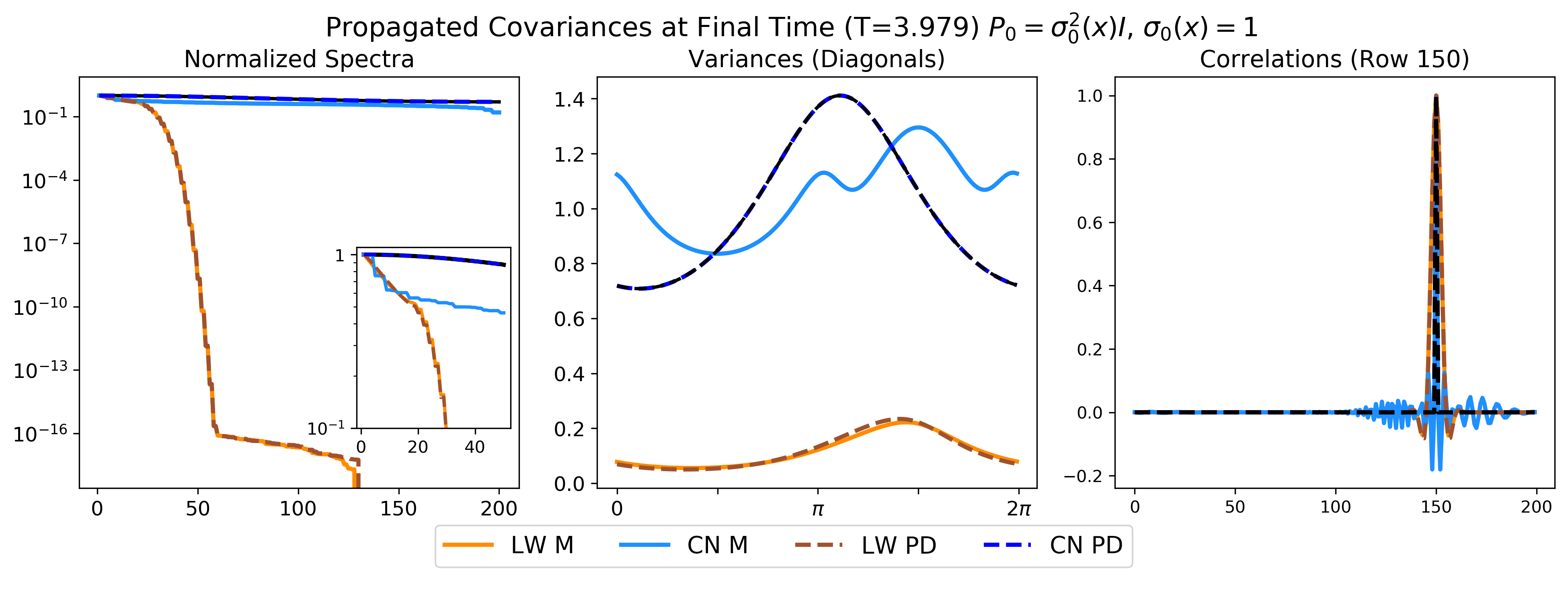}
    \label{fig:tri ID b}\includegraphics[width=\linewidth]{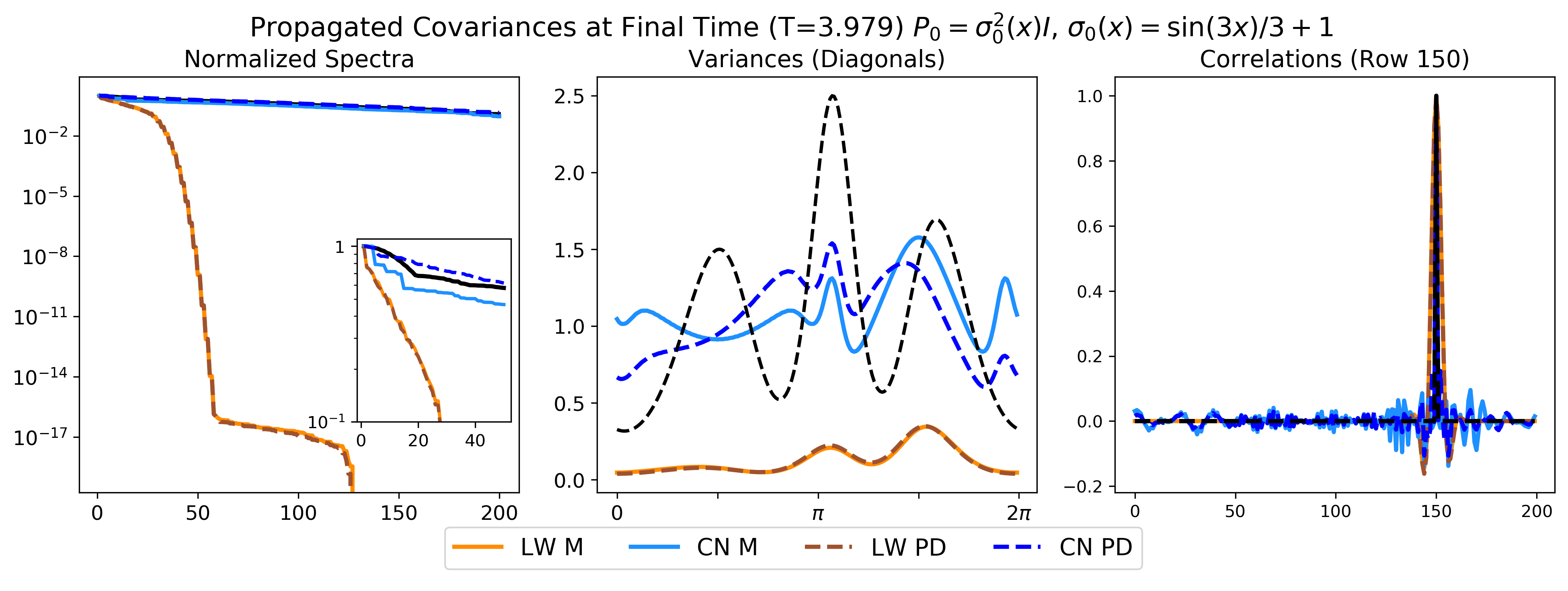}
    \caption{Same as \cref{fig:GC0.25} for covariances with zero initial correlation length ($c = L = 0$).}
    \label{fig:tri ID}
\end{figure}

In the limiting case when the initial correlation lengths become zero, we see two contrasting behaviors in the numerically propagated covariances depending on the initial variance; see \cref{fig:tri ID}. When the initial variance is constant (i.e. the initial covariance is the identity matrix), the Crank-Nicolson polar decomposition is the only scheme that correctly captures the behavior of the exact covariance (top row of \cref{fig:tri ID}). This is expected from the continuum analysis coupled with the quadratically conservative property of the Crank-Nicolson scheme used to construct the matrix $\bs{U}$. When the initial variance varies spatially (bottom row of \cref{fig:tri ID b}), the Crank-Nicolson polar decomposition propagation instead behaves more similarly to the Crank-Nicolson traditional propagation and both have regions of variance loss and variance gain.  
Carefully comparing the traditional and polar decomposition diagonals, the Crank-Nicolson polar decomposition propagation is a slight improvement over the Crank-Nicolson traditional propagation, but these differences are relatively minor compared to their absolute errors. The Lax-Wendroff schemes are substantially dissipative in all cases. We also observe that as the values of $c$ and $L$ decrease towards zero, the normalized spectra in \cref{fig:GC1,fig:GC0.25,fig:FOAR0.25,fig:tri ID} decay more slowly and become relatively flat. This suggests that low-rank approximations would have difficulty capturing these covariances as correlation lengths shrink. 

\begin{figure}[htbp]
    \centering
    \includegraphics[width=\linewidth]{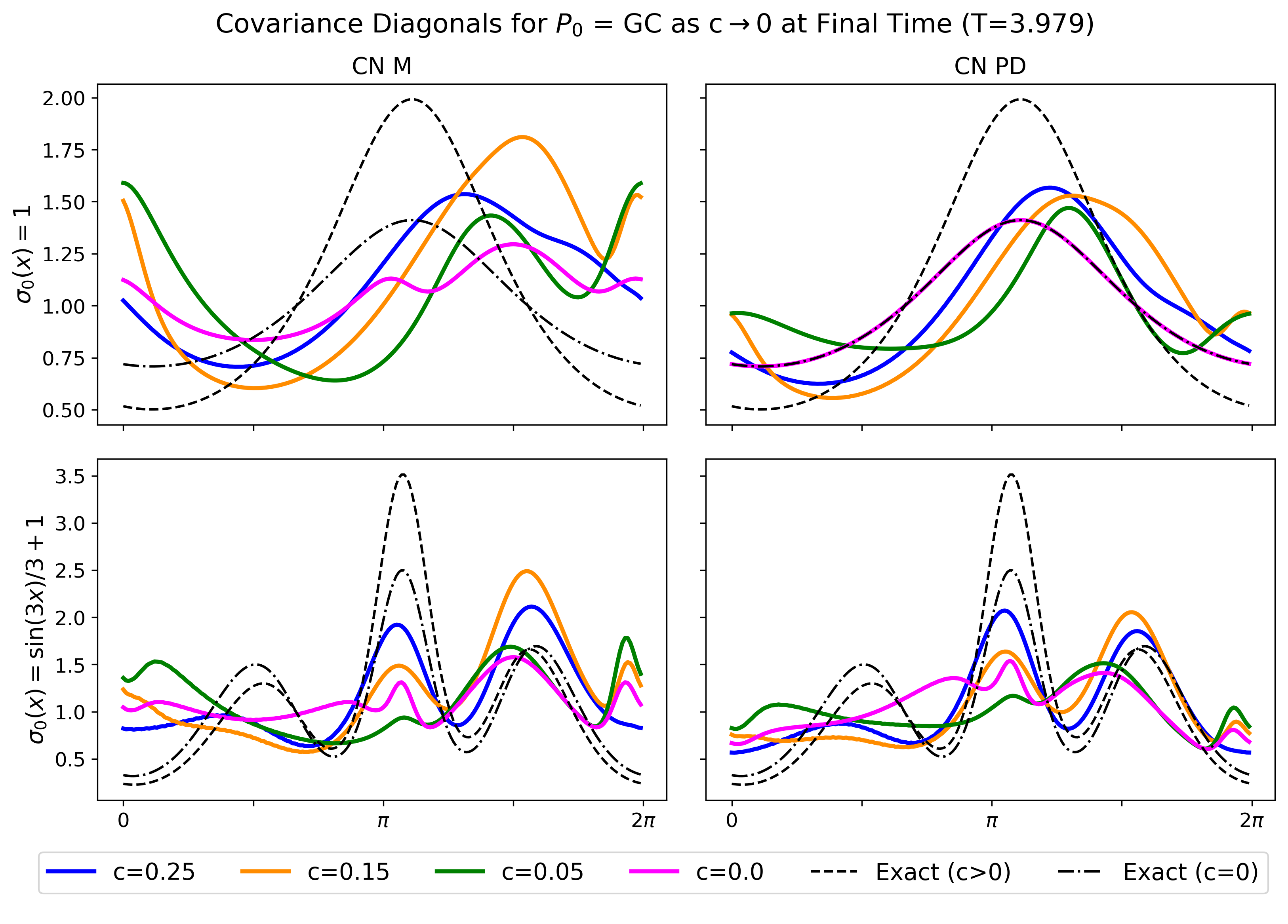}
    \caption{Covariance diagonals extracted from Crank-Nicolson traditional and polar decomposition propagation methods at the final time $T$ for GC initial correlations as $c$ approaches zero. Top row: constant initial variance. Bottom row: spatially-varying initial variance. Left column: propagation via traditional Crank-Nicolson. Right column: propagation via Crank-Nicolson polar decomposition. Black curves are the exact diagonals, dashed for covariances with nonzero initial correlation lengths (variance solution) and dot-dashed for covariances with zero initial correlation length (continuous spectrum solution). These exact curves here labeled as {\rm Exact} $(c>0)$ and {\rm Exact} $(c=0)$ are the same exact curves labeled as {\rm Exact} $\sigma^2(x,t)$ and {\rm Exact} $P^c(x,t)$ in \cref{fig:variance and cts spec}.
    In the top right panel, the exact curve for $c=0$ (black dot-dashed) and CN PD curve for $c=0$ (magenta) identically overlap.}
    \label{fig:GC diags}
\end{figure}

\begin{figure}[htbp]
    \centering
    \includegraphics[width=\linewidth]{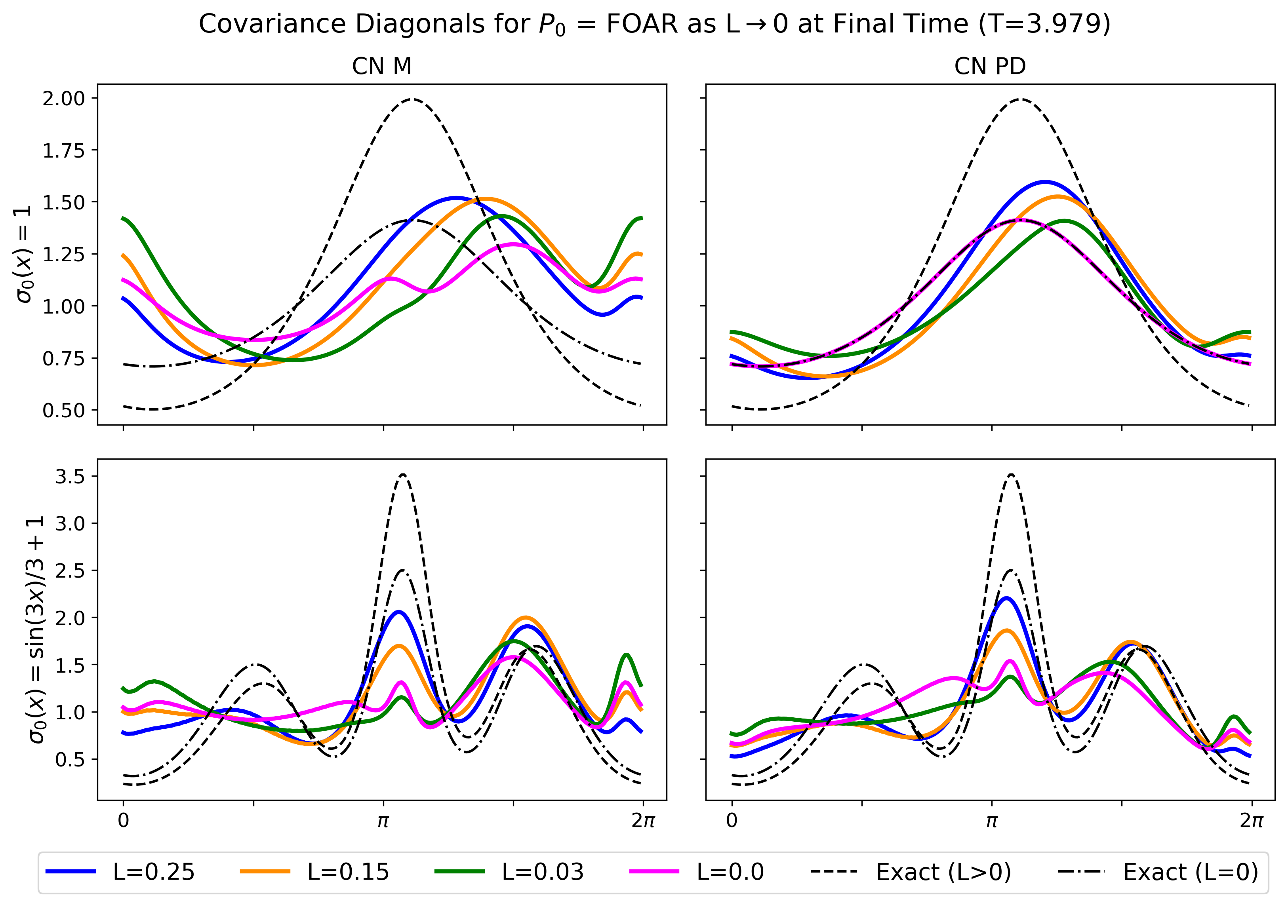}
    \caption{Same as \cref{fig:GC diags} for FOAR correlation function initial covariances.}
    \label{fig:FOAR diags}
\end{figure}

Comparing the diagonals at the final time across a series of initial correlation lengths in \cref{fig:GC diags,fig:FOAR diags} demonstrates the severity of the variance loss and gain caused by inaccurate variance propagation and provides a closer look at the approach to a limiting behavior seen in the trace time series. Without prior knowledge of the discontinuous change in diagonal dynamics as the initial correlation length tends to zero, one might surmise from \cref{fig:GC diags,fig:FOAR diags} that the observed diagonal behavior is caused by dissipation. Knowing the continuum behavior, however, makes it clear that we are observing inaccurate variance propagation associated with the discontinuous change in continuum dynamics. The Crank-Nicolson polar decomposition propagation for a constant initial variance is the only scheme that captures the correct diagonal behavior when the initial covariance is the identity, and gradually approaches this behavior as the initial correlation length decreases. For all other cases in \cref{fig:GC diags,fig:FOAR diags}, the resulting variances are smooth solutions, but incorrect ones, and gradually approach their own limiting behavior at $c=L=0$ rather than changing abruptly as in the continuum case. The Lax-Wendroff schemes are severely dissipative, hardly resembling the correct diagonal dynamics, and are not shown in these and subsequent figures.

\begin{figure}[htbp]
    \centering
    \includegraphics[width=\linewidth]{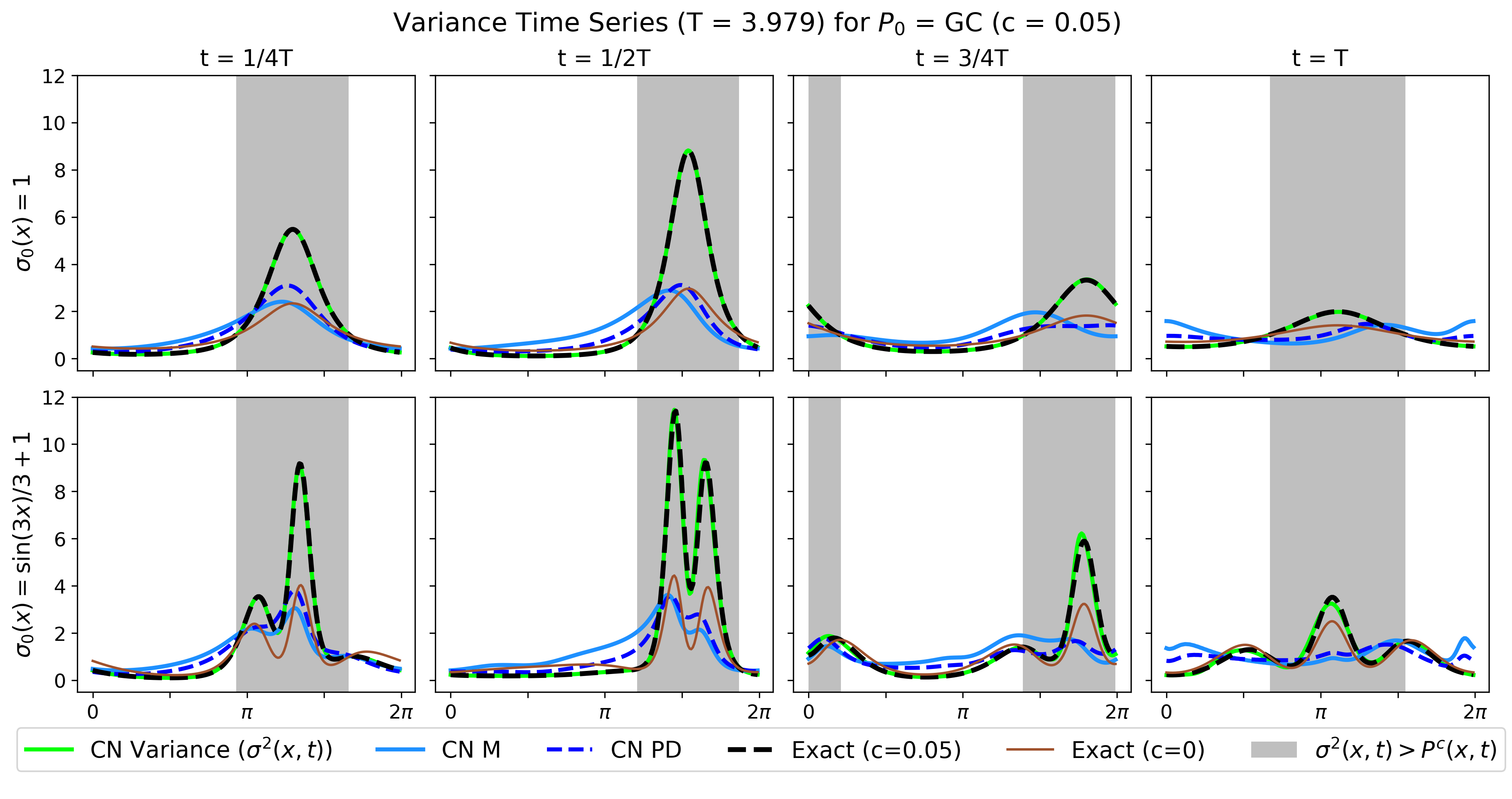}
    \caption{Variance time series for GC case with $c=0.05$. Top row: constant initial variance. Bottom row: spatially-varying initial variance. Crank-Nicolson traditional and polar decomposition propagation are shown (solid and dashed blue) as well as propagating the variance by solving the variance equation (one-dimensional version of \cref{eq:variance S2}) independently using Crank-Nicolson (green). The exact solutions to the variance and continuous spectrum equations  are shown in black-dashed and solid brown and denoted as  {\normalfont Exact} $(c=0.05)$ and {\normalfont Exact} $(c=0)$, respectively. Regions where the exact variance is larger (smaller) than the exact continuous spectrum are highlighted in grey (white).}
    \label{fig:GC var}
\end{figure}

The errors in the diagonal propagation are not limited to the final time; errors start to accumulate early on in the propagation cycle. We show this for the GC case in \cref{fig:GC var}, where the FOAR results are similar but not shown. We can see clearly that rather than approximating the variance for $c>0$ (black dashed), the Crank-Nicolson schemes (blue) tend to approximate the diagonal for $c=0$ (brown), e.g. the continuous spectrum, which is not the correct diagonal behavior for this case. Along with the variances extracted from the propagated covariances in \cref{fig:GC var}, we include the diagonal propagated independently by solving the one-dimensional version of \cref{eq:variance S2} using the Crank-Nicolson scheme, as was shown in \cref{fig:variance and cts spec}. Including the variance solution computed using the Crank-Nicolson scheme in \cref{fig:GC var} emphasizes that propagating the covariance diagonal independently using the known diagonal dynamics significantly reduces the errors in the variance compared to the variance extracted from the propagated covariance matrix. 

Grey regions in \cref{fig:GC var} correspond to where the exact variance solution is larger than the exact continuous spectrum solution, and these regions tend to correspond to where the numerically propagated diagonal for the short, nonzero initial correlation length (solid light blue) exhibits variance loss. Further insight into this behavior can be gained from the continuum perspective as follows. Assuming $P_0^d(\bs{x},\bs{x}) = P^c_0(\bs{x})$, as is the case in \cref{fig:GC var}, one can show, using the generalized variance equation \cref{eq:variance} and continuous spectrum equation \cref{eq:continuous spectrum}, that the two satisfy the following relation,
\begin{equation}\label{eq:m relation}
    \frac{\sigma^2 (\bs{x},t)}{P^c (\bs{x},t)} = m(\bs{x},t),
\end{equation}
where $m = m(\bs{x},t)$ satisfies the continuity equation with unit initial condition,
\begin{gather}
    m_t + \bs{v\cdot\nabla}m + (\bs{\nabla\cdot v})m = 0, \nonumber \\
    m(\bs{x},t_0) = 1. \label{eq:m continuity}
\end{gather}
If we consider \cref{eq:m continuity} in the context of mass conservation, the ratio $\sigma^2/P^c$ must be conserved. Therefore in regions of mass convergence ($m >1$) we have $\sigma^2 > P^c$, and conversely in regions of mass divergence ($m<1$) we have $\sigma^2 < P^c$. In fact, the function $m$ of \cref{eq:m continuity} on $S^1_1$ with advection speed \cref{eq:1d advection speed} can be expressed explicitly using \cref{eq:1d continuity exact soln} taking one as the initial condition. Solving $m = 1$, or equivalently $v(s(x,t)) - v(x) = 0$ (following the notation of \cref{sec:appendix b}), determines the boundaries between regions of mass convergence and divergence in $S^1_1$ at time $t > 0$.

For short initial correlation lengths, if the numerical schemes are better approximating the continuous spectrum solution $P^c$, regions of mass convergence in the ratio $\sigma^2/P^c$ should correspond to variance loss and regions of mass divergence in the ratio $\sigma^2/P^c$ should correspond to variance gain. Regions where the variance (black-dashed) is larger than the continuous spectrum (brown), e.g. mass convergence, highlighted in grey in \cref{fig:GC var}, generally do coincide with regions where we see variance loss in the Crank-Nicolson traditional and polar decomposition schemes. Conversely, in the unshaded regions of \cref{fig:GC var} where the variance is less than the continuous spectrum, e.g. mass divergence, the Crank-Nicolson traditional and polar decomposition schemes generally exhibit variance gain. Hence, we can exploit the ratio in \cref{eq:m relation} to quantify variance loss and gain in relation to mass convergence and divergence.

Since the function $m$ must satisfy mass conservation, there will always be regions of mass convergence and mass divergence in the spatial domain at a given time $t$. In the context of advective systems, this implies that we should see both loss and gain of variance over time as a consequence of the ratio in \cref{eq:m relation}; whether there will be global variance loss or variance gain depends on the advection speed $\bs{v}$. In all the numerical results, \cref{fig:GC1,fig:trace trig,fig:GC0.25,fig:FOAR0.25,fig:tri ID,fig:GC diags,fig:FOAR diags,fig:GC var}, we see both spurious variance loss and gain as a result of inaccurate variance propagation. Data assimilation literature tends to focus on loss of variance because of its negative impact on data assimilation schemes \cite[section 8.8]{Ja70}, \cite[section 9.2]{Ma82}, \cite{AnAn99,HoMi00,LyCoZhChMeOlRe04}. However, we see here for a general advective system that both loss and gain exist because of conservation, cumulatively causing inaccurate variance propagation. Variance inflation, a tool often used in data assimilation practice to combat variance loss by rescaling the variance, typically by a multiplicative or additive factor, could be optimally adjusted in this context. Distinguishing between regions of mass convergence and divergence in general isolates the regions of variance loss and gain, thus the inflation factor can be tuned to only inflate in regions of mass convergence where we expect variance loss. This would avoid using a single scale-factor that inflates the whole variance function and prevent variance inflation in regions where we already see variance gain.

Since covariance information on and close to the diagonal may be sufficient information for many applications, local covariance evolution, where the variance and correlation lengths are propagated rather than the full matrix, may prove useful. \cite{Co93} first discussed local covariance evolution through continuum analysis of hyperbolic and parabolic PDEs, similar to the equations discussed here. The covariance equation in $2N$ spatial dimensions, where $N$ is the number of spatial dimensions of the state, is reduced to a system of auxiliary PDEs in $N$ dimensions consisting of variance and correlation length equations, which  approximate the covariance locally. Our results, together with the results of \cite{Co93}, suggest further investigation into local covariance propagation, which may help reduce computational expense and the spurious loss and gain of variance observed in full covariance propagation.

\section{Conclusions and discussion}
\label{sec:conclusions and discussion}
In this work, we study covariance propagation associated with random state variables governed by a generalized advection equation. We do this in an effort to understand the root causes of spurious loss of variance, which is observed in data assimilation schemes wherein the covariance is explicitly or implicitly propagated. Using a continuum analysis to guide the interpretation of our numerical results, new insights are gained through the detailed study of both the continuum and discrete covariance evolution.  The main conclusions are as follows:

\begin{enumerate}
    \item Continuum analysis of the state and covariance equations is necessary to establish a fundamental understanding of the covariance evolution. In particular, the continuum analysis uncovers the discontinuous behavior of the diagonal dynamics as the correlation length approaches zero, for example in the vicinity of sharp gradients, which is an insight crucial for understanding the spurious loss and gain of variance observed in our numerical experiments.
   
    \item Comparison of the numerical results with the continuum analysis shows that full-rank covariance propagation via \cref{eq:discrete covariance} typically results in considerable spurious loss of variance. This is due to the peculiar discontinuous behavior of the dynamics far more than to any numerical dissipation. Discrete propagation using \cref{eq:discrete covariance} produces inaccurate variances, resulting in both variance loss in regions of mass convergence and variance gain in regions of mass divergence. When propagating initial covariances with short, nonzero correlation lengths, the numerical schemes better approximate the dynamics of the zero correlation length case  than those of the nonzero correlation length case.
    
    
    \item Isolating the variance and propagating it independently ameliorates the variance loss and gain observed during full-rank covariance propagation and yields accurate propagation because it adheres to the continuum dynamics. This result suggests further investigation into alternate methods of covariance propagation. 
    
\end{enumerate}

Though the continuum analysis performed here may not be tractable in all situations, it serves as a foundation for understanding covariance evolution and interpreting the results of our numerical experiments. Central to the continuum analysis is the polar decomposition of the fundamental solution operator \cref{eq:m polar decomp}, which is general in that it is a canonical form for all bounded linear operators on Hilbert spaces. Thus, we expect that much of the work presented here can be extended to general hyperbolic systems of PDEs having a quadratic energy functional, e.g. \cite{Co10}. 

It is important to recognize that the loss of variance observed in our numerical experiments is a result of the discontinuous change in continuum covariance dynamics discussed in \cref{sec:continuum}. Even when propagating covariance matrices using a fully Lagrangian scheme, as done in \cite{LyCoZhChMeOlRe04}, propagated covariances still suffer from spurious loss of variance that is not due to the numerical scheme but rather to the discontinuous change in dynamics that we've identified in this work. Covariance propagation tends to be overlooked in the data assimilation literature as a potential source of variance loss, particularly when using the same numerical method that propagates the state. Data assimilation schemes that do not propagate the covariance explicitly may experience errors similar to what we observe here because the underlying cause of these errors is the  covariance dynamics, not the numerical scheme. Studying full-rank covariance propagation as in \cref{eq:discrete covariance} isolates the spurious loss of variance as an issue with covariance dynamics and implies that approximations to \cref{eq:discrete covariance}, such as in ensemble Kalman filters, can suffer variance loss in a similar manner. The errors caused during the covariance propagation may be a neglected source of the model error or ``system error" observed in data assimilation schemes \cite[p.\,3285]{HoMi05}. 

Our work brings to light a fundamental issue associated with current approaches to numerical covariance propagation and recommends investigation into alternate methods of covariance propagation. \Cref{fig:variance and cts spec,fig:GC var}, which display the independent variance propagation for example, suggest that local covariance evolution may be an adequate alternative. As discussed in \cite{Co93}, for applications in which information on and close to the diagonal is sufficient, evolving the variance and correlation lengths themselves may serve as an alternative to full covariance evolution. 

\appendix
\section{Proofs of solutions with Dirac delta initial conditions}\label{sec:appendix a}
We begin by verifying that $P(\bs{x}_1,\bs{x}_2,t) = d^2(\bs{x}_1,t)\delta(\bs{x}_1,\bs{x}_2)$, where $d^2$ satisfies \cref{eq:d2 equation}, is the solution to the covariance evolution equation \cref{eq:continuum covariance} for $P_0(\bs{x}_1,\bs{x}_2) = \delta(\bs{x}_1,\bs{x}_2)$.
\begin{proof}
    The function $P(\bs{x}_1,\bs{x}_2,t) = d^2(\bs{x}_1,t)\delta(\bs{x}_1,\bs{x}_2)$ is a weak (distribution) solution of \cref{eq:continuum covariance} with $P_0(\bs{x}_1,\bs{x}_2) = \delta(\bs{x}_1,\bs{x}_2)$ if
    \begin{equation}\label{eq:weak soln def}
        (\mathcal{L}^*\phi,P) = 0
    \end{equation}
    for all test functions $\phi \in C^1(\Omega\times\Omega\times [t_0,T])$ with period $t=T$, where $\mathcal{L}^*$ is the adjoint differential operator
    \begin{equation}\label{eq:adjoint differential operator}
        \mathcal{L}^* = -\partial_t - \bs{v}_1\bs{\cdot\nabla}_1- \bs{v}_2\bs{\cdot\nabla}_2 + b_1 + b_2 - \bs{\nabla}_1\bs{\cdot v}_1 - \bs{\nabla}_2\bs{\cdot v}_2
    \end{equation} corresponding to the differential operator
    \begin{equation}\label{eq:differential operator}
        \mathcal{L} = \partial_t + \bs{v}_1\bs{\cdot\nabla}_1 + \bs{v}_2\bs{\cdot\nabla}_2 + b_1 + b_2
    \end{equation}
    of \cref{eq:continuum covariance}. Here we denote $(\cdot,\cdot)$ as the inner product over $L^2(\Omega\times\Omega\times [t_0,T])$ and we will denote $(\cdot,\cdot)'$ as the inner product over $L^2(\Omega\times [t_0,T])$. 
    Substituting $P = d^2(\bs{x}_1,t)\delta(\bs{x}_1,\bs{x}_2)$ into the expression for $(\mathcal{L}^*\phi,P)$, applying the Dirac delta, and expanding the result yields
    \begin{equation}\label{eq:soln line 2}
        (\mathcal{L}^*\phi, P) =  (-\phi_t(\bs{x}_1,\bs{x}_1,t) - \bs{\nabla}_1\bs{\cdot}(\bs{v}_1\phi(\bs{x}_1,\bs{x}_1,t)) + (2b_1-\bs{\nabla}_1\bs{\cdot v}_1)\phi(\bs{x}_1,\bs{x}_1,t),d^2(\bs{x}_1,t))',
    \end{equation}
    where we note that $\bs{\nabla}_2\bs{\cdot v}_2|_{\bs{x}_2 = \bs{x}_1} = \bs{\nabla}_1\bs{\cdot v}_1$ and $\bs{\nabla}_1\phi(\bs{x}_1,\bs{x}_1,t) = \bs{\nabla}_1\phi(\bs{x}_1,\bs{x}_2,t)|_{\bs{x}_2 = \bs{x}_1} + \bs{\nabla}_2\phi(\bs{x}_1,\bs{x}_2,t)|_{\bs{x}_2 = \bs{x}_1}$. Using integration by parts to move derivatives off of the test function $\phi$ onto $d^2$ in \cref{eq:soln line 2} gives us
    \begin{equation}\label{eq:soln line 3}
        (\mathcal{L}^*\phi, P)= (\phi(\bs{x}_1,\bs{x}_1,t),d^2_t + \bs{v}_1\bs{\cdot\nabla}_1d^2 + (2b_1-\bs{\nabla}_1\bs{\cdot v}_1)d^2)' = 0,
    \end{equation}
    since $d^2$ satisfies \cref{eq:d2 equation}. Thus,  $P(\bs{x}_1,\bs{x}_2,t) = d^2(\bs{x}_1,t)\delta(\bs{x}_1,\bs{x}_2)$ is a weak solution to \cref{eq:continuum covariance}.
\end{proof}

To show that $\tilde{P} = \tilde{P}^c(\bs{x}_1,t)\delta(\bs{x}_1,\bs{x}_2)$ is the solution to \cref{eq:P tilde equation} for $\tilde{P}_0 = P_0^c(\bs{x}_1)\delta(\bs{x}_1,\bs{x}_2)$ follows an identical proof as given above, where we take
\begin{align}
    b_i &\mapsto \frac12 \bs{\nabla}_i\bs{\cdot v}_i, \ i=1,2 \\
    d^2(\bs{x}_1,t) &\mapsto \tilde{P}^c(\bs{x}_1,t),
\end{align}
where $\tilde{P}^c$ satisfies \cref{eq:P tilde c}.

\section{Solution to the state equation}\label{sec:appendix b}
In this appendix we solve the one-dimensional version of \cref{eq:continuity S2} used in the numerical experiments with advection speed \cref{eq:1d advection speed}, 
\begin{gather}
    q_t + vq_x + v'q = 0, \ x \in S_1^1,\nonumber \\
    q(x,t_0) = q_0(x), \label{eq:1d continuity} \\
    v = \sin(x)+2. \nonumber
\end{gather}

To solve \cref{eq:1d continuity}, we start by rewriting the equation in terms of the Lagrangian/total derivative
\begin{equation}\label{eq:lagrangian derivative}
    \frac{D}{Dt} \equiv \partial_t + v\partial_x.
\end{equation}
By doing so, \cref{eq:1d continuity} becomes the ordinary differential equation,
\begin{gather}
    \frac{Dq(x(t),t)}{Dt} = -v'(x(t))q(x(t),t), \nonumber \\
    q(x,t_0) = q_0(x). \label{eq:lagrangian continuity}
\end{gather}
Associated with \cref{eq:lagrangian continuity} are the Lagrangian trajectories, or characteristic equations, that determine how the spatial variable $x$ evolves over time,
\begin{gather}
    \frac{dx(t)}{dt} = v(x(t)), \nonumber \\
    x(t_0) = s. \label{eq:characteristic equation}
\end{gather}
Here, we take $s$ as a general initial parameter and define our solutions to \cref{eq:characteristic equation} as $x(t;s)$, following the notation of \cite[chapter 3]{ShLe15}. 

Since the state equation is reduced to an ordinary differential equation in \cref{eq:lagrangian continuity}, solutions are of the form
\begin{equation}\label{eq:continuity soln 1}
    q(x,t) = q_0(x(t_0))e^{\int_{t_0}^t -v'(x(\tau))d\tau}.
\end{equation}
We can solve the integral in the exponential of \cref{eq:continuity soln 1} using a simple $u-$substitution of $u = x(t)$ along with \cref{eq:characteristic equation},
\begin{equation}\label{eq:exponential integral}
    \int_{t_0}^t -v'(x(\tau))d\tau =\int_{s}^{x(t)} -\frac{v'(u)}{v(u)}du = \ln\left(\frac{v(s)}{v(x(t))}\right).
\end{equation}

To rewrite \cref{eq:exponential integral} in terms of $(x,t)$, we first solve for $x(t;s)$ using the characteristic equation \cref{eq:characteristic equation}, then invert to determine $s = s(x,t)$, which can be done since the advection speed $v$ is continuously differentiable \cite[chapter 3]{ShLe15}. Therefore, substituting \cref{eq:exponential integral} and $x(t_0) = s(x,t)$ into \cref{eq:continuity soln 1}, we have the explicit solution to \cref{eq:1d continuity},
\begin{gather}
    q(x,t) = q_0(s(x,t))\left(\frac{\sin(s(x,t))+2}{\sin(x)+2}\right), \label{eq:1d continuity exact soln}\\
    s(x,t) = 2\tan^{-1}\left(\frac{\sqrt{3}}{2}\tan\left(\tan^{-1}\left(\frac{2\tan(x/2)+1}{\sqrt{3}}\right) - \frac{\sqrt{3}t}{2}\right) - \frac12\right).\label{eq:s}
\end{gather}
The function $\tan^{-1}$ refers to the principle branch of inverse tangent, where $y=\tan^{-1}(x),\ x \in \mathbb{R}, \ y \in (-\pi/2,\pi/2)$, defined by solving $x=\tan(y)$ \cite[section 6.3.3]{Zw96}. 

\begin{figure}[htbp]
    \centering
    \includegraphics[width=\linewidth]{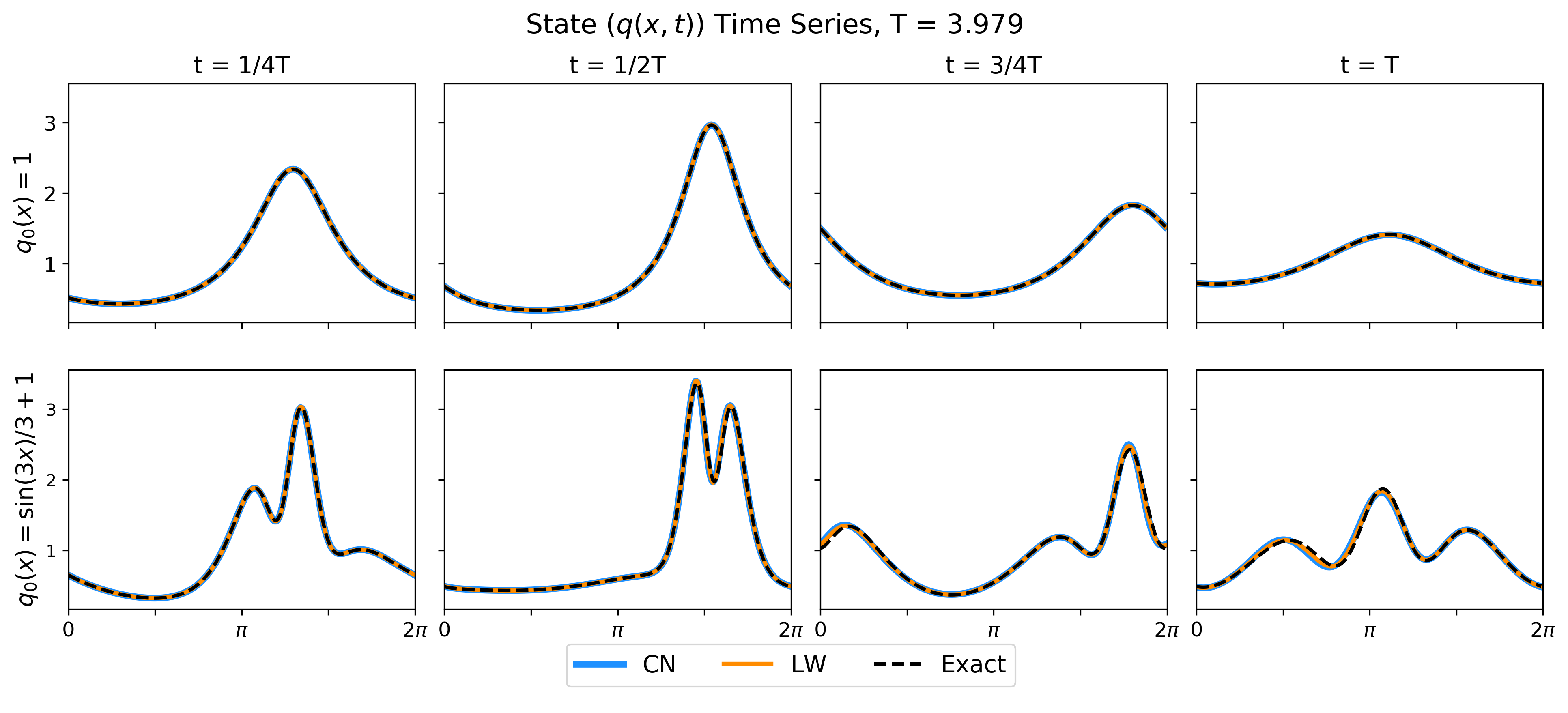}
    \caption{Solutions to \cref{eq:1d continuity} for $q_0(x) = 1$ (top row) and $q_0(x) = \sin(3x)/3 + 1$ (bottom row) at four different times. The exact solution (black dashed) is shown along with finite-difference solutions using Crank-Nicolson (blue) and Lax-Wendroff (orange). The finite difference solutions are nearly indistinguishable from the exact solution.}
    \label{fig:state prop}
\end{figure}

\Cref{fig:state prop} plots \cref{eq:1d continuity exact soln} for two different initial conditions, along with solving \cref{eq:1d continuity} using Lax-Wendroff and Crank-Nicolson schemes on a uniform grid of 200 grid points and unit Courant-Friedrichs-Lewy number.

The function $m$ of \cref{eq:m continuity} on $S^1_1$ with advection speed \cref{eq:1d advection speed} can be expressed explicitly using \cref{eq:1d continuity exact soln} taking one as the initial condition. The exact solution $d$ of \cref{eq:d equation} on $S_1^1$ for $b = v_x$, as in the numerical experiments, can also be computed explicitly from \cref{eq:1d continuity exact soln} (as $d^2$ in this case satisfies \cref{eq:1d continuity} with unit initial condition). 

\section*{Acknowledgments}
This material is based upon work supported by the National Science Foundation (NSF) Graduate Research Fellowship Program under Grant No.~DGE-1650115. TM is supported by the NSF CAREER Program under Grant No.~AGS-1848544. SEC is supported by the National Aeronautics and Space Administration (NASA) under Modeling, Analysis and Prediction (MAP) program Core funding to the Global Modeling and Assimilation Office (GMAO) at Goddard Space Flight Center (GSFC). Any opinions, findings, and conclusions or recommendations expressed in this material are those of the author(s) and do not necessarily reflect the views of the NSF or of NASA. 
\bibliographystyle{siamplain}
\bibliography{references}

\end{document}

%% file: Gilpin_covprop_juq_shared.tex
\usepackage{lipsum}
\usepackage{amsfonts}
\usepackage{graphicx}
\usepackage{epstopdf}
\usepackage{algorithmic}
\ifpdf
  \DeclareGraphicsExtensions{.eps,.pdf,.png,.jpg}
\else
  \DeclareGraphicsExtensions{.eps}
\fi

\usepackage{enumitem}
\setlist[enumerate]{leftmargin=.5in}
\setlist[itemize]{leftmargin=.5in}

\newsiamremark{remark}{Remark}
\newsiamremark{hypothesis}{Hypothesis}
\crefname{hypothesis}{Hypothesis}{Hypotheses}
\newsiamthm{claim}{Claim}

\headers{Continuum Covariance Propagation for Understanding Variance Loss in Advective Systems}{S. Gilpin, T. Matsuo, and S. E. Cohn}

\title{Continuum Covariance Propagation for Understanding Variance Loss in Advective Systems\thanks{Submitted to the editors August 25, 2021.
\funding{This work was funded by the National Science Foundation under grant No.~DGE-1650115 and AGS-1848544, and by the National Aeronautics and Space Administration (NASA) Modeling, Analysis and Prediction (MAP) program under Core funding to the Global Modeling and Assimilation Office (GMAO) at Goddard Space Flight Center (GSFC).}}}

\author{Shay Gilpin\thanks{Department of Applied Mathematics, University of Colorado Boulder, Boulder, CO, Ann and H.J. Smead Department of Aerospace Engineering Sciences, University of Colorado Boulder, Boulder, CO, (\email{shay.gilpin@colorado.edu}).}
\and Tomoko Matsuo\thanks{Ann and H.J. Smead Department of Aerospace Engineering Sciences, University of Colorado Boulder, Boulder, CO, Department of Applied Mathematics, University of Colorado Boulder, Boulder, CO, (\email{tomoko.matsuo@colorado.edu}).}
\and Stephen E. Cohn\thanks{Global Modeling and Assimilation Office, NASA Goddard Space Flight Center, Greenbelt, MD}, (\email{stephen.e.cohn@nasa.gov})}

\usepackage{amsopn}
